\documentclass[10pt]{article}
\usepackage[english]{babel}

\usepackage{latexsym}
\usepackage{amsmath}
\usepackage{amsfonts}
\usepackage{mathrsfs}
\usepackage{amstext}
\usepackage{amssymb}
\usepackage{amsthm}       %proof environment
\usepackage{pgf,tikz}
\usetikzlibrary{arrows}

\usepackage{wrapfig}

\usepackage{xspace}

\newtheorem{thm}{Theorem}[section]
\newtheorem{prop}[thm]{Proposition}

\newtheorem{defn}[thm]{Definition}

\newtheorem{lem}[thm]{Lemma}

\usepackage[margin=3cm]{geometry}

\begin{document}

\title{Existence results for a nonlinear transmission problem}
\author{M.~Dalla Riva \footnote{Centre for Research and Development in Mathematics and Applications (CIDMA) - University of Aveiro - Portugal } \quad \& \quad G.~Mishuris \footnote{Department of Mathematics - Aberystwyth University - UK}}
\maketitle

\begin{abstract}
Let $\Omega^o$ and $\Omega^i$ be open bounded subsets of $\mathbb{R}^n$ of class $C^{1,\alpha}$ such that the closure  of $\Omega^i$ is contained in $\Omega^o$. Let $f^o$ be a function in $C^{1,\alpha}(\partial\Omega^o)$ and let $F$ and $G$  be continuous functions from $\partial\Omega^i\times\mathbb{R}$ to $\mathbb{R}$.  By exploiting an argument based on potential theory and on the Leray-Schauder principle we show that under suitable and completely explicit conditions on $F$ and $G$ there exists at least one pair of continuous functions $(u^o, u^i)$ such that
\[
\left\{
\begin{array}{ll}
\Delta u^o=0&\text{in }\Omega^o\setminus\mathrm{cl}\Omega^i\,,\\
\Delta u^i=0&\text{in }\Omega^i\,,\\
u^o(x)=f^o(x)&\text{for all }x\in\partial\Omega^o\,,\\
u^o(x)=F(x,u^i(x))&\text{for all }x\in\partial\Omega^i\,,\\
\nu_{\Omega^i}\cdot\nabla u^o(x)-\nu_{\Omega^i}\cdot\nabla u^i(x)=G(x,u^i(x))&\text{for all }x\in\partial\Omega^i\,,
\end{array}
\right.
\]
 where the last equality is attained in certain weak sense. In  a simple example we show that such a pair of functions $(u^o, u^i)$ is in general neither unique nor local unique. If instead the fourth condition of the problem is obtained by a small nonlinear perturbation of a homogeneous linear condition, then we can prove the existence of at least one classical solution which is in addition locally unique.  \\

\medskip

\end{abstract}

\section{Introduction}
We investigate the existence of solutions for a boundary value problem with a nonlinear transmission condition. In order to define such a boundary value problem we  introduce some notation. We fix once for all   
\[
\text{a natural number $n\in\mathbb{N}$, $n\ge 2$, and a real number $\alpha\in]0,1[$,}
\] 
where $\mathbb{N}$ denotes the set of natural numbers including $0$. Then we fix two sets $\Omega^o$ and $\Omega^i$ in the $n$-dimensional Euclidean space $\mathbb{R}^n$. The   letter `$o$' stands for `outer domain' and the letter `$i$' stands for `inner domain'. We assume that $\Omega^o$ and $\Omega^i$ satisfy the following condition:
\[
\begin{split}
&\text{$\Omega^o$ and $\Omega^i$ are open bounded subsets of $\mathbb{R}^n$ of class $C^{1,\alpha}$, $\mathrm{cl}\Omega^i\subseteq\Omega^o$, and}\\
&\text{the boundaries $\partial\Omega^o$ and $\partial\Omega^i$ are connected.}
\end{split}
\]
For the definition of functions and sets of the usual Schauder class $C^{0,\alpha}$ and $C^{1,\alpha}$, we refer for example to Gilbarg and Trudinger \cite[\S6.2]{GiTr01}. Here and in the sequel $\mathrm{cl}\Omega$ denotes the closure of $\Omega$ for all $\Omega\subseteq\mathbb{R}^n$. 
Then we fix a function $f^o\in C^{1,\alpha}(\partial\Omega^o)$ and two  continuous  functions $F$ and $G$ from $\partial\Omega^i\times \mathbb{R}$ to $\mathbb{R}$ and consider the following nonlinear transmission boundary value problem for a pair of functions $(u^o,u^i)$  in $C^{1,\alpha}(\mathrm{cl}\Omega^o\setminus\Omega^i)\times C^{1,\alpha}(\mathrm{cl}\Omega^i)$,
\begin{equation}\label{nlprob}
\left\{
\begin{array}{ll}
\Delta u^o=0&\text{in }\Omega^o\setminus\mathrm{cl}\Omega^i\,,\\
\Delta u^i=0&\text{in }\Omega^i\,,\\
u^o(x)=f^o(x)&\text{for all }x\in\partial\Omega^o\,,\\
u^o(x)=F(x,u^i(x))&\text{for all }x\in\partial\Omega^i\,,\\
\nu_{\Omega^i}\cdot\nabla u^o(x)-\nu_{\Omega^i}\cdot\nabla u^i(x)=G(x,u^i(x))&\text{for all }x\in\partial\Omega^i\,,
\end{array}
\right.
\end{equation}
where $\nu_{\Omega^i}$ denotes the outer unit normal to the boundary $\partial{\Omega^i}$. Our  aim is to determine suitably general and completely explicit conditions on  $F$ and $G$ which ensure the existence of solutions of problem \eqref{nlprob}.

The analysis of problems such as \eqref{nlprob} is motivated by the role played in continuum mechanics. In particular, nonlinear transmission conditions of this kind arise in the study of composite structures  glued together by thin adhesive layers which are thermically or mechanically very different from the components' constituents. In modern material technology such composites are widely used (see, {\it e.g.}, \cite{MiMiOc08, MiMiOc09, RoCa01}), but  the numerical treatment of the mathematical model by finite elements methods is still difficult, requires the introduction of highly inhomogeneous meshes, and often leads to poor accuracy and numerical instability (see, {e.g.}, Babu\v{s}ka and Suri \cite{BaSu92}). A convenient way to overcome this problem is to replace the thin layers  by zero thickness interfaces between the composite's components. Then  one  has to define on such interfaces suitable transmission conditions which incorporates the thermical and mechanical properties of the original layers.  Such a procedure can be rigorously justified by an asymptotic method   and leads to the introduction of boundary value problems with nonlinear transmission conditions such as those in \eqref{nlprob} (see for example \cite{MiOc13} and the references therein).  

We observe that the existence of solutions of nonlinear boundary value problems has been largely investigated by means of variational techniques (see, {\it e.g.}, the monographs of Ne\v{c}as \cite{Ne83} and of Roub\'i\v{c}ek \cite{Ro13} and the references therein). In fact, under some restrictive assumptions on the functions $F$ and $G$, the existence of solutions of our problem \eqref{nlprob} could be  deduced by exploiting some known results. In particular, if it happens that problem \eqref{nlprob} can be reformulated into an   equation of the form $-\mathrm{div} A(x, U)\nabla U= 0$, where $A$ is a suitable Carath\'eodory function and the unknown function $U$ belongs to the Sobolev space $H^1(\Omega^o)$ and satisfies a Dirichlet condition on $\partial
\Omega^o$,  then the existence and uniqueness of a solution   can be directly deduced by the results of Hlav\'a\v{v}cek, K\v{r}\'i\v{z}ek and Mal\'y in \cite{HlKrMa94}. This is for example the case when $G=0$ and the function $F(x,t)$ of $(x,t)\in \partial\Omega^i\times\mathbb{R}$  is constant with respect to $x$, is differentiable with respect to $t$, and the partial differential $\partial_t F(x_0,\cdot)$ is Lipschitz continuous and satisfies the inequality $1/c<\partial_t F(x_0,t)<c$ for a constant $c>0$ and for all $t\in\mathbb{R}$ (here $x_0$ is a fixed point of $\partial\Omega^i$). 

In this paper instead, we exploit a method based on potential theory to rewrite problem  \eqref{nlprob} into a suitable nonlinear system of integral equations which can be analysed by a fixed-point theorem.  Potential theoretic techniques have been largely exploited in literature to study  existence and uniqueness problems for linear or semilinear partial differential equations with non linear boundary conditions. In particular, as far back as in 1921 Carleman \cite{Ca21} has considered the existence of harmonic functions $u$ in a domain $\Omega$ which satisfy  a non-linear Robin condition $\nu_{\Omega}(x)\cdot\nabla u(x)=H(x,u(x))$ on the boundary $\partial\Omega$.  Since then, such a problem has received the attention of many authors such as Leray \cite{Le33} (see also Jacob \cite{Ja35}), Nakamori and Suyama \cite{NaSu50},  Kilngelh\"ofer \cite{Kl68, Kl69}, Cushing \cite{Cu71}, and Efendiev, Schmitz, and Wendland \cite{EfScWe99}. In the case of domains with a small hole we also  mention the  nonlinear Robin problem for the Laplace operator investigated in Lanza de Cristoforis \cite{La07} and the nonlinear traction problem in elasticity addressed in \cite{DaLa10}.   Moreover, an approach based on coupling of boundary integral and finite element methods has been developed in order to study  exterior nonlinear boundary value problems with  transmission conditions, we mention for example the papers of Berger \cite{Be89},  Berger, Warnecke, and Wendland \cite{BeWaWe90}, Costabel and Stephan \cite{CoSt90},  and Gatica and Hsiao \cite{GaHs92, GaHs95}. In particular, Barrenechea and Gatica considered in \cite{BaGa96} the case when the jump of the normal derivative across the interface boundary depends nonlinearly on the Dirichlet data. Boundary integral methods have been applied also by Mityushev and Rogosin for the analysis of transmission problems in the two dimensional plane (cf.~\cite[Chap.~5]{MiRo00}). Finally, we mention the nonlinear transmission problem in a domain with a small inclusion investigated by  Lanza de Cristoforis in \cite{La10} and the periodic analog studied by Lanza de Cristoforis and Musolino in \cite{LaMu14}.  

\section{Description of the main results}

  We now describe the main results of the present paper. We will exploit the following notation: if $H$ is a function from $\partial\Omega^i\times\mathbb{R}$, then we denote by $\mathcal{F}_H$   the nonlinear non-autonomous composition operators which takes a function $f$ from $\partial\Omega^{i}$ to $\mathbb{R}$ to the function $\mathcal{F}_H f$  defined by 
\[
\mathcal{F}_Hf(x)\equiv H(x,f(x))\qquad\forall x\in\partial\Omega^i\,.
\]
Since the functions $F$ and $G$ which define the nonlinear condition in \eqref{nlprob} are assumed to be continuos from $\partial\Omega^i\times\mathbb{R}$ to $\mathbb{R}$, one easily verifies that $\mathcal{F}_F$ and $\mathcal{F}_G$ are continuous from $C^0(\partial\Omega^i)$ to itself. Then we consider the following condition: 
\begin{equation}\label{assumption0}
\text{the composition operator $(I_{\Omega^i}+\mathcal{F}_F)$ has a continuous inverse $(I_{\Omega^i}+\mathcal{F}_F)^{(-1)}$ from $C^0(\partial\Omega^i)$ to itself.}
\end{equation}
Here $I_{\Omega^i}$ denotes the identity operator from $C^0(\partial\Omega^i)$ to itself. 
We observe that for the validity of condition \eqref{assumption0} it is not required that the  function which takes $t$ to $F(x,t)$ is monotone for all fixed $x\in\partial\Omega^i$. 
In addition, we introduce a condition on the magnitude of $F$ and $G$: we assume that 
\begin{equation}\label{assumption1}
\begin{split}
&\text{there exist $c_1, c_2\in ]0,+\infty[$,  $\delta_1\in ]1,+\infty[$, and $\delta_2\in [0,1[$  such that}\\
&|F(x,t)|\ge c_1|t|^{\delta_1}-(1/c_1)\quad\text{ and  }\\
 & |G(x,t)|\le c_2(1+|F(x,t)|)^{\delta_2}\qquad\forall(x,t)\in\partial\Omega^i\times\mathbb{R}\,.
\end{split}
\end{equation} The first condition in \eqref{assumption1} is a super-linear grow condition for $F$, while the second one is a sub-linear  grow condition for $G$ with respect to $F$  (which  is a strictly weaker condition than the standard sub-linear condition $|G(x,t)|\le c_2(1+|t|)^{\delta_2}$). 

By exploiting an argument based on the invariance of the Leray-Schauder topological degree we show in our main Theorem \ref{existence}   that   conditions \eqref{assumption0} and \eqref{assumption1} imply the  existence of at least one pair of continuous functions $(\tilde u^o, \tilde u^i)\in C^0(\mathrm{cl}\Omega^o\setminus\Omega^i)\times C^0(\mathrm{cl}\Omega^o)$ which satisfies the first four equations of \eqref{nlprob} in the classical sense and fulfils the fifth condition in a certain weak sense which will be clarified (see Definition \ref{weak} below). However,  the conditions in \eqref{assumption0} and \eqref{assumption1} do not imply neither the uniqueness nor the local uniqueness of the pair $(\tilde u^o, \tilde u^i)$. 

This last fact can be evidenced in a simple example. Take $\Omega^o=R\mathbb{B}_n$,  $\Omega^i=r\mathbb{B}_n$, with $r,R\in\mathbb{R}$, $r<R$, and  $\mathbb{B}_n\equiv\{x\in\mathbb{R}^n\,:\;|x|<1\}$. Then assume that $f^o$ is constant and identically equal to a real number  $t^o\in\mathbb{R}$ and that $F(x,t)\equiv f(t)$ and $G(x,t)\equiv g(t)$ for all $(x,t)\in\partial\Omega^i\times\mathbb{R}$, where $f$ and $g$ are continuous functions from $\mathbb{R}$ to $\mathbb{R}$. We set 
\begin{figure}%{h}{0.4\textwidth}
\begin{center}
\includegraphics[width=5cm]{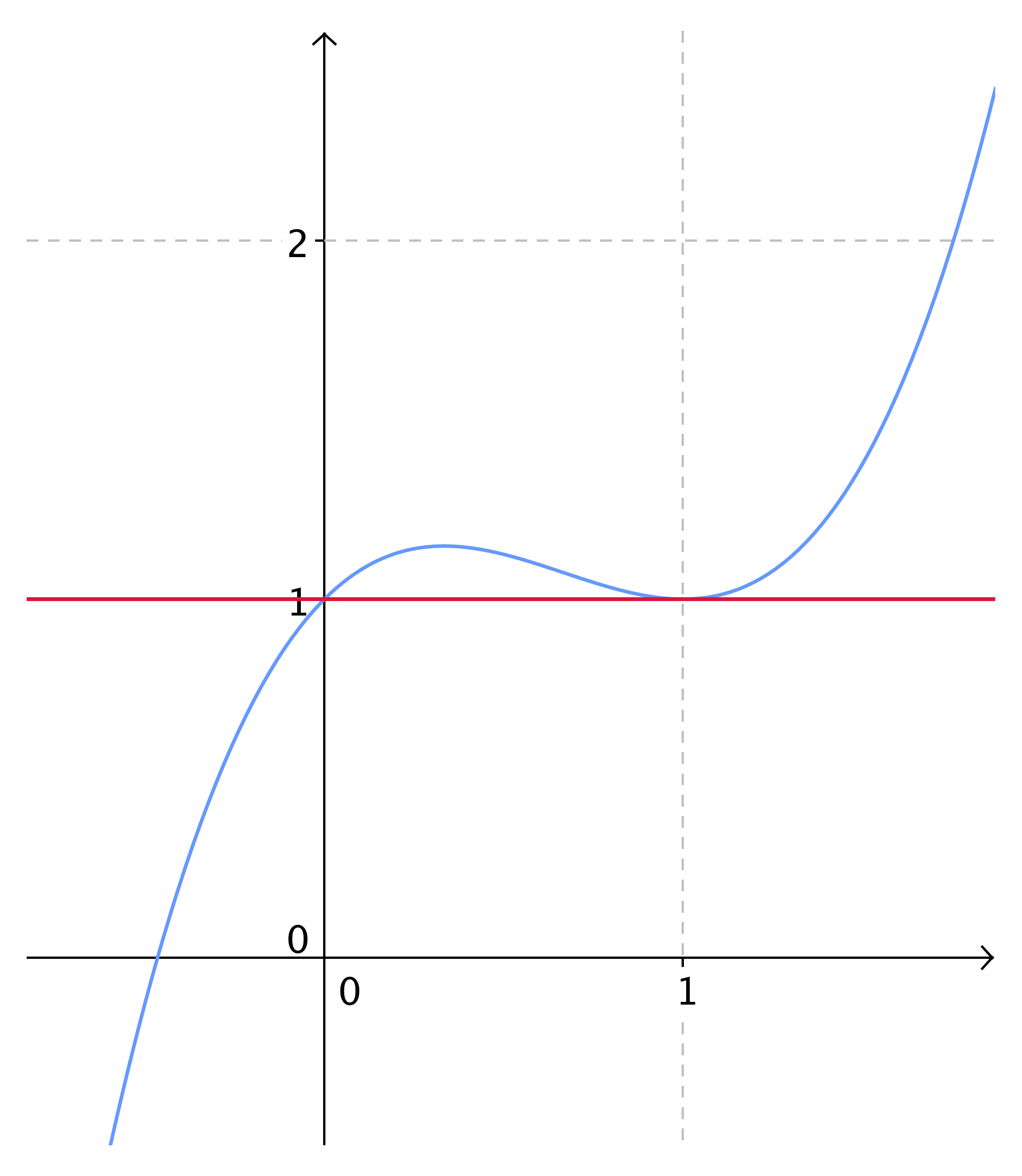}
\caption{the intersections of the blue graph with the red line correspond to solutions of \eqref{1d}}
\label{fig1}
\end{center}
\end{figure}
\[
\Gamma_n(t)\equiv
\left\{
\begin{array}{ll}
\frac{1}{2\pi}\log t&\text{if }n=2\,,\\
\frac{t^{2-n}}{s_n(2-n)}&\text{if }n\ge 3\,,
\end{array}
\right. \quad\forall t\in]0,+\infty[\,,
\] 
where $s_n$ denotes the $(n-1)$-dimensional measure of $\partial\mathbb{B}_n$ (thus $\Gamma_n(|x|)=S_n(x)$ with  $S_n$ the standard fundamental solution of $\Delta$, see also definition \eqref{Sn} below). Then the pair of functions $(u^o,u^i)$ defined by 
\begin{equation}\label{uiuo1d}
\begin{split}
&u^o(x)=t^o-\frac{\Gamma_n(R)-\Gamma_n(|x|)}{\Gamma'_n(r)}g(t^i)\quad\forall x\in\mathrm{cl}\Omega^o\setminus\Omega^i\,,\\
&u^i(x)=t^i\quad\forall x\in\mathrm{cl}\Omega^i
\end{split}
\end{equation}
is a solution of problem \eqref{nlprob} for all $t^i\in\mathbb{R}$ which are solutions of the equation
 \begin{equation}\label{1d}
t^o-\frac{\Gamma_n(R)-\Gamma_n(r)}{\Gamma'_n(r)}g(t^i)=f(t^i)\,.
 \end{equation}
Now take
\[
f(t)\equiv t^3-2t^2+t+1\qquad\forall t\in\mathbb{R}
\] 
and assume that $g$ is constant. One immediately verifies that the corresponding functions $F$ and $G$ satisfy the conditions in \eqref{assumption0} and \eqref{assumption1}. In addition, if $t^o$, $R$, $r$, and $g$ are choosen in such a way that the left hand side of \eqref{1d} is equal to $1$,  then equation \eqref{1d} has two solutions: $t^i=0$ and $t^i=1$ (see Fig.\ref{fig1}). Accordingly,  the corresponding problem \eqref{nlprob} has at least two different solutions provided by \eqref{uiuo1d}. If instead  $f(t)\equiv t^3-2t^2+t+1$ for $t<0$ and $t>1$ and $f(t)\equiv 1$ for $t\in[0,1]$, then every $t^i$ in $[0,1]$ is a solution of \eqref{1d} and the corresponding  solutions of problem \eqref{nlprob} are not locally unique in any reasonable topology.

We observe that Theorem \ref{existence} shows the existence of pair of functions $(\tilde u^o,\tilde u^i)\in C^0(\mathrm{cl}\Omega^o\setminus\Omega^i)\times C^0(\mathrm{cl}\Omega^i)$ which are solutions of problem \eqref{nlprob} in a certain `weak' sense but it would be preferable to have classical solutions in  
$C^{1,\alpha}(\mathrm{cl}\Omega^o\setminus\Omega^i)\times C^{1,\alpha}(\mathrm{cl}\Omega^i)$ (or at least in $H^1(\Omega^o\setminus\mathrm{cl}\Omega^i)\times H^1(\Omega^i)$).  Thus, it is natural to ask what further conditions should one impose on $F$ and $G$ in order to obtain such a regularity. In Theorem \ref{existenceC0a} we show that, if
\begin{equation}\label{assumption2}
\text{$(I_{\Omega^i}+\mathcal{F}_F)^{(-1)}$ and $\mathcal{F}_G$ map $C^{0,\alpha}(\partial\Omega^i)$ to itself,}
\end{equation} then problem \eqref{nlprob} has at least one weak solution in $C^{0,\alpha}(\mathrm{cl}\Omega^o\setminus\Omega^i)\times C^{0,\alpha}(\mathrm{cl}\Omega^o)$. However, in order to obtain solutions in  $C^{1,\alpha}(\mathrm{cl}\Omega^o\setminus\Omega^i)\times C^{1,\alpha}(\mathrm{cl}\Omega^i)$ by exploiting our argument it does not suffice to increase the regularity of $F$ and $G$ and it seems that a different approach should be implemented. 

To illustrate this fact, we consider in the last Section \ref{sec.small} the case when the fourth condition of problem \eqref{nlprob} is a small nonlinear perturbation of a homogenous linear condition.  Namely,  we assume that $F(x,t)=\lambda t+\epsilon\Phi(x,t)$ for all $(x,t)\in\partial\Omega^i\times\mathbb{R}$, where $\lambda$ is a positive real constant, $\epsilon$ is a small real parameter, and $\Phi$ is a continuous function from $\partial\Omega^i\times\mathbb{R}$ to $\mathbb{R}$. Then we consider the nonlinear transmission problem
\begin{equation}\label{nlprob.small}
\left\{
\begin{array}{ll}
\Delta u^o=0&\text{in }\Omega^o\setminus\mathrm{cl}\Omega^i\,,\\
\Delta u^i=0&\text{in }\Omega^i\,,\\
u^o(x)=f^o(x)&\text{for all }x\in\partial\Omega^o\,,\\
u^o(x)=\lambda u^i(x)+\epsilon   \Phi(x,u^i(x))&\text{for all }x\in\partial\Omega^i\,,\\
\nu_{\Omega^i}\cdot\nabla u^o(x)-\nu_{\Omega^i}\cdot\nabla u^i(x)=G(x,u^i(x))&\text{for all }x\in\partial\Omega^i\,,\\
\end{array}
\right.
\end{equation} 
\begin{figure}%{h}%{0.4\textwidth}
\begin{center}
\includegraphics[width=6cm]{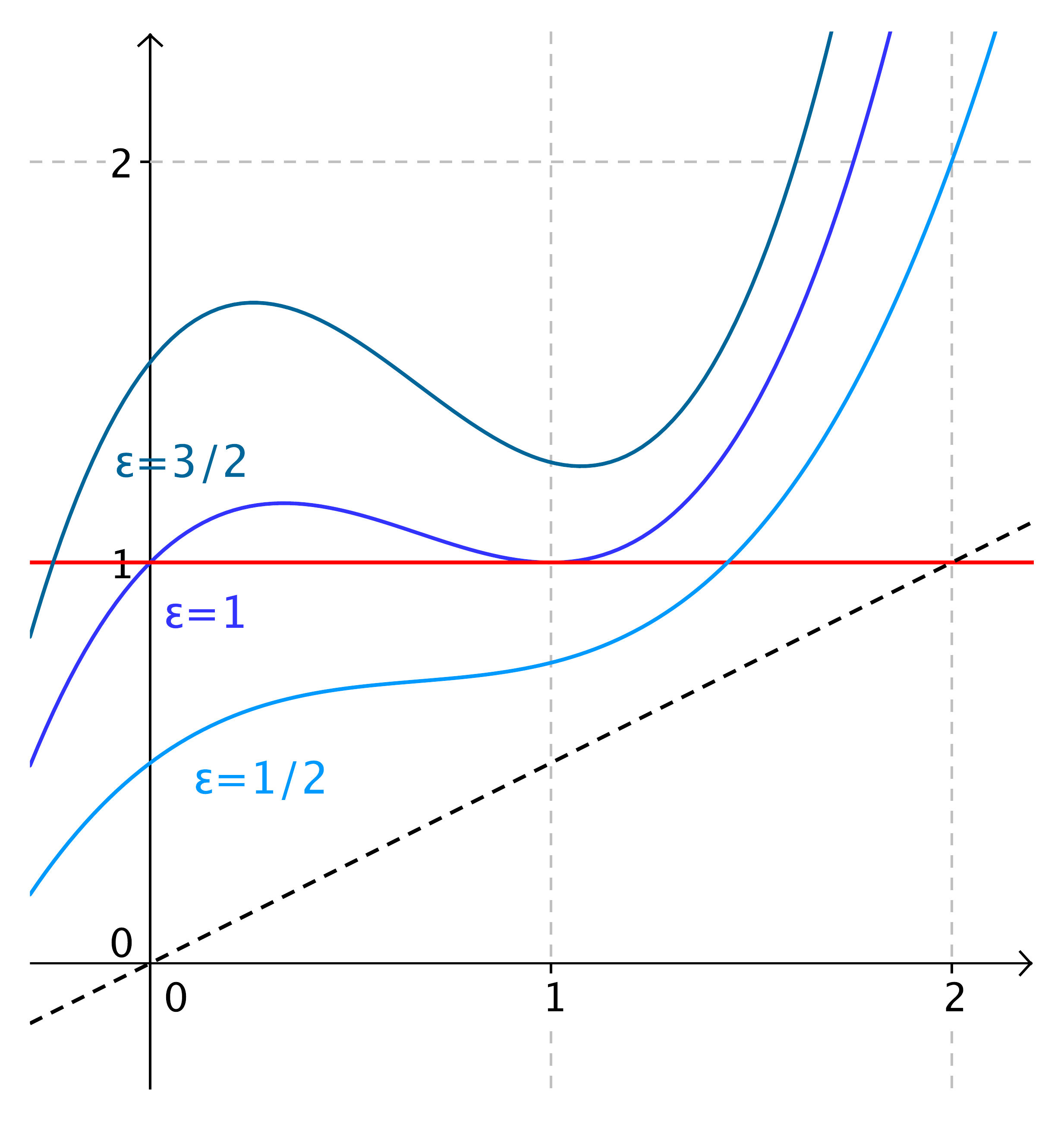}
\caption{the blue lines are graphs of $\lambda t+\epsilon \phi(t)$ for $\epsilon=3/2$, $\epsilon=1$, and $\epsilon=1/2$.   The intersections of the blue graphs with the red line correspond to solutions of \eqref{1dsmall}}
\label{fig2}
\end{center}
\end{figure} 
for a pair of functions $(u^o,u^i)\in C^{1,\alpha}(\mathrm{cl}\Omega^o\setminus\Omega^i)\times C^{1,\alpha}(\mathrm{cl}\Omega^i)$.  In Theorem \ref{thm.small} below we show that, under suitable assumptions on $\Phi$, $G$, and $\Omega^i$ (see condition \eqref{ass.omega}), there exists $\epsilon_*>0$ such that problem \eqref{nlprob.small} has a solution $(u^o_\epsilon,u^i_\epsilon)\in C^{1,\alpha}(\mathrm{cl}\Omega^o\setminus\Omega^i)\times C^{1,\alpha}(\mathrm{cl}\Omega^i)$ for all $\epsilon\in]-\epsilon_*,\epsilon_*[$. Such a solution $(u^o_\epsilon,u^i_\epsilon)$ is locally unique in $C^{1,\alpha}(\mathrm{cl}\Omega^o\setminus\Omega^i)\times C^{1,\alpha}(\mathrm{cl}\Omega^i)$ for all fixed $\epsilon\in]-\epsilon_*,\epsilon_*[$ and, in addition, the map which takes $\epsilon$ to $(u^o_\epsilon,u^i_\epsilon)$ is continuously Fr\'echet differentiable from $]-\epsilon_*,\epsilon_*[$ to $C^{1,\alpha}(\mathrm{cl}\Omega^o\setminus\Omega^i)\times C^{1,\alpha}(\mathrm{cl}\Omega^i)$.   However, Theorem \ref{thm.small} does not provide any estimate for the value of $\epsilon_*$. Therefore, the existence conditions provided by Theorem \ref{thm.small} are not completely explicit, as instead are those of Theorem \ref{existence}. 

We also observe  that   the assumptions of Theorems \ref{existence} and  \ref{thm.small} may be simultaneously verified, but the solutions provided by   Theorem \ref{existence} may not coincide with those provided by Theorem \ref{thm.small}. Consider for example the case introduced here above where $\Omega^o=R\mathbb{B}_n$ and $\Omega^i=r\mathbb{B}_n$, with $r<R$. Assume that $f^o(x)=t^o$, $\Phi(x,t)=\phi(t)$, and $G(x,t)=g(t)$ for all $(x,t)\in\partial\Omega^i\times\mathbb{R}$, where $t^o\in\mathbb{R}$ and $\phi$ and $g$ are continuous functions from $\mathbb{R}$ to itself. Then we look for solutions of problem \eqref{nlprob.small} in the form \eqref{uiuo1d} with $t^i\in\mathbb{R}$ solution of the equation  
 \begin{equation}\label{1dsmall}
t^o-\frac{\Gamma_n(R)-\Gamma_n(r)}{\Gamma'_n(r)}g(t^i)=\lambda t+\epsilon\phi(t^i)\,.
 \end{equation}
Now we take $\lambda\equiv 1/2$, 
\[
\phi(t)\equiv t^3-2t^2+\frac{1}{2}t+1\qquad\forall t \in\mathbb{R}\,,
\]
and $g$ constant. One can choose $t^o$, $r$, $R$, and $g$ in such a way that the left hand side of \eqref{1dsmall} is equal to $1$. Then it is easily verified that equation \eqref{1dsmall} has two solutions for $\epsilon=1$: $t^i=0$ and $t^i=1$. Instead, for $\epsilon>1$ we only have the solution provided by  Theorem \ref{existence} and due to the behaviour at infinity of $\phi$  and we loose the solution provided by Theorem \ref{thm.small} and due to the smallness of $\epsilon$ (see Fig.~\ref{fig2}).   Similar examples can be prepared to show that the local uniqueness of the solution guaranteed by Theorem \ref{thm.small} for $\epsilon$ small can be lost for $\epsilon=1$.

Finally, we observe that potential theoretic methods have been developed by Escauriaza {\it et al.}~\cite{EsFaVe92, EsSe93, EsMi04} for the analysis of linear transmission problems  in Lipschitz domains. However, the argument used in the present paper for the proof of  the main Theorem \ref{existence}  cannot be immediately extended to the case of a  Lipschitz contact boundary $\partial\Omega^i$. The reason is that  the compactness of the double layer operator $W_{\Omega^i}$ plays a crucial role in the proof of Proposition \ref{existenceC}, where we apply the Leray-Schauder principle to prove that the fixed point equation \eqref{integral} has solutions (see Section \ref{sec.preliminaries} for the definition of $W_{\Omega^i}$). As is well known,   $W_{\Omega^i}$ is compact in $L^p(\partial\Omega^i)$,  $p\in]1,+\infty[$, if $\Omega^i$ is at least of class $C^1$, but may be not compact if $\Omega^i$ is just a Lipschitz domain (cf., {\it e.g.}, Fabes {\it et al.}~\cite{FaJoLe77, FaJoRi78}).  

The paper is organised as follows. Section \ref{sec.preliminaries} is a section of preliminaries where we introduce some classical notion of potential theory. In Section \ref{sec.existence} we prove our main Theorem \ref{existence} where we show the existence of continuous solutions of problem \eqref{nlprob}. Finally, in Section \ref{sec.small} we consider problem \eqref{nlprob.small} and we show the existence of locally unique $C^{1,\alpha}$ solutions for $\epsilon$ small.

\section{Classical notions of potential theory}\label{sec.preliminaries}

We denote by $S_n$ the function from $\mathbb{R}^n\setminus\{0\}$ to $\mathbb{R}$ defined by
\begin{equation}\label{Sn}
S_n(x)\equiv
\left\{
\begin{array}{ll}
\frac{1}{2\pi}\log\,|x|&\text{if $n=2$,}\\
\frac{1}{s_n(2-n)}|x|^{2-n}&\text{if $n\ge 3$,}
\end{array}
\right.
\qquad\forall x\in\mathbb{R}^2\setminus\{0\}\,.
\end{equation}
 As is well known, $S_n$ is a fundamental solution for the Laplace operator in $\mathbb{R}^n$.

Let $\Omega$ be an open bounded subset of $\mathbb{R}^n$ of class $C^{1,\alpha}$. Let $\phi\in L^2(\partial\Omega)$.  Then $v_\Omega[\phi]$ denotes the single layer potential with density $\phi$. Namely,
\[
v_\Omega[\phi](x)\equiv\int_{\partial\Omega}\phi(y)S_n(x-y)\,d\sigma_y\qquad\forall x\in\mathbb{R}^n\\,
\] where $d\sigma$ denotes the area  element on $\partial\Omega$. As is well known, if $\phi\in L^\infty(\partial\Omega)$, then  $v_\Omega[\phi]$ is a continuous function from $\mathbb{R}^n$ to $\mathbb{R}$. In addition, if  $\phi\in C^{0,\alpha}(\partial\Omega)$, then the restrictions $v^+_\Omega[\phi]\equiv v_\Omega[\phi]_{|\mathrm{cl}\Omega}$ and $v^-_\Omega[\phi]\equiv v_\Omega[\phi]_{|\mathbb{R}^n\setminus\Omega}$ belong to $C^{1,\alpha}(\mathrm{cl}\Omega)$ and to $C^{1,\alpha}_{\mathrm{loc}}(\mathbb{R}^n\setminus\Omega)$, respectively.  Here $C^{1,\alpha}_{\mathrm{loc}}(\mathbb{R}^n\setminus\Omega)$ denotes the space of functions on $\mathbb{R}^n\setminus\Omega$ whose restrictions to $\mathrm{cl}\mathcal{O}$ belong to  $C^{1,\alpha}(\mathrm{cl}\mathcal{O})$ for all open bounded subsets $\mathcal{O}$ of $\mathbb{R}^n\setminus\Omega$.

If  $\psi\in L^2(\partial\Omega)$, then $w_\Omega[\psi]$ denotes the double layer potential with density $\psi$. Namely,
\[
w_\Omega[\psi](x)\equiv-\int_{\partial\Omega}\psi(y)\;\nu_{\Omega}(y)\cdot\nabla S_n(x-y)\,d\sigma_y\qquad\forall x\in\mathbb{R}^n\,,
\] where $\nu_\Omega$ denotes the outer unit normal to $\partial\Omega$  and the symbol `$\cdot$' denotes the scalar product in $\mathbb{R}^n$. If $\psi\in C^{1,\alpha}(\partial\Omega)$, then the restriction $w_\Omega[\psi]_{|\Omega}$ extends to a function $w^+_\Omega[\psi]$ of $C^{1,\alpha}(\mathrm{cl}\Omega)$ and  the restriction $w_\Omega[\psi]_{|\mathbb{R}^n\setminus\mathrm{cl}\Omega}$ extends to a function $w^-_\Omega[\psi]$ of $C^{1,\alpha}_{\mathrm{loc}}(\mathbb{R}^n\setminus\Omega)$.

 Let 
\[
W_\Omega[\psi](x)\equiv -\int_{\partial\Omega}\psi(y)\;\nu_{\Omega}(y)\cdot\nabla S_n(x-y)\, d\sigma_y\qquad\forall x\in\partial\Omega\,,
\] for all $\psi\in L^2(\partial\Omega)$, and
\[
W^*_\Omega[\phi](x)\equiv \int_{\partial\Omega}\phi(y)\;\nu_{\Omega}(x)\cdot\nabla S_n(x-y)\, d\sigma_y\qquad\forall x\in\partial\Omega\,,
\] for all $\phi\in L^2(\partial\Omega)$.  As is well known $W_\Omega$ and $W^*_\Omega$ are compact operator  from  $L^2(\partial\Omega)$ to itself and are adjoint one to the other. In the sequel we denote by  $I_\Omega$  the identity map from $L^2(\partial\Omega)$ to itself. Thus 
\[
\text{$\pm\frac{1}{2}I_\Omega+ W_\Omega$ and $\pm\frac{1}{2}I_\Omega+ W^*_\Omega$ are Fredholm operators of index $0$ from $L^2(\partial\Omega)$ to itself.}
\]

We now introduce  the following classical result of Schauder \cite{Sc31, Sc32}:
\begin{lem}\label{schauder}
Let $\beta\in]0,1]$. Then the map which takes $\psi$ to $W_\Omega[\psi]$ is continuous from $C^0(\partial\Omega)$ to  $C^{0,\alpha}(\partial\Omega)$ and from $C^{1,\beta}(\partial\Omega)$ to  $C^{1,\alpha}(\partial\Omega)$. The map which takes $\phi$ to $W^*_\Omega[\phi]$ is continuous from $C^{0,\beta}(\partial\Omega)$ to  $C^{0,\alpha}(\partial\Omega)$.
\end{lem}
As a consequence, the map which takes $\psi$ to $W_\Omega[\psi]$ is compact from   $C^{1,\alpha}(\partial\Omega)$  to itself and the map which takes $\phi$ to $W^*_\Omega[\phi]$  is compact from  $C^{0,\alpha}(\partial\Omega)$ to itself. Then one immediately deduces the validity of the following.
\begin{lem}\label{fredC}
The operators $\pm\frac{1}{2}I_\Omega+ W_\Omega$ are Fredholm of index $0$ from $C^0(\partial\Omega)$ to itself, from $C^{0,\alpha}(\partial\Omega)$ to itself, and from $C^{1,\alpha}(\partial\Omega)$ to itself. The operators $\pm\frac{1}{2}I_\Omega+ W^*_\Omega$ are Fredholm of index $0$ from  $C^{0,\alpha}(\partial\Omega)$ to itself.
\end{lem}
In addition we have the following technical Lemma \ref{regularity}.
\begin{lem}\label{regularity}
Let $\psi\in L^2(\partial\Omega)$. Let $\beta\in[0,\alpha]$.  Let $\gamma\in\mathbb{R}$.  If   $(\frac{1}{2}I_\Omega+\gamma W^*_{\Omega})\psi$ belongs to $C^{0,\beta}(\partial\Omega)$, then $\psi\in C^{0,\beta}(\partial\Omega)$.  
\end{lem}
\proof If $(\frac{1}{2}I_\Omega+\gamma W^*_{\Omega})\psi\in C^{0,\beta}(\partial\Omega)$, then a standard argument based on iterated kernels ensures that $\psi\in C^0(\partial\Omega)$. It follows that $W^*_{\Omega}\psi\in C^{0,\beta'}(\partial\Omega)$ for all $\beta'\in[0,\alpha[$ (cf.~Miranda \cite[Chap.~II, \S14, IV]{Mi70}, see also Schauder \cite{Sc32}). Thus $\psi=2(\frac{1}{2}I_\Omega+\gamma W^*_{\Omega})\psi-2\gamma W^*_{\Omega}\psi$ belongs to  $C^{0,\beta''}(\partial\Omega)$ with $\beta''\equiv\min\{\beta',\beta\}$ for all $\beta'\in[0,\alpha[$. Accordingly the lemma is proved for $\beta<\alpha$. If instead $\beta=\alpha$, then we observe that $W^*_{\Omega}\psi\in C^{0,\alpha}(\partial\Omega)$ by the membership of $\psi$ in $C^{0,\beta'}(\partial\Omega)$ with $\beta'\in]0,1]$ (cf.~Lemma \ref{schauder}). Thus $\psi=2(\frac{1}{2}I_\Omega+\gamma W^*_{\Omega})\psi-2\gamma W^*_{\Omega}\psi$ belongs to  $C^{0,\alpha}(\partial\Omega)$  and the Lemma is proved.\qed 

\vspace{\baselineskip}

A similar result holds if we replace in Lemma \ref{regularity} the operator $W^*_{\Omega}$ by the operator $W_{\Omega}$ (cf.~Miranda \cite[Chap.~II, \S15, II]{Mi70}). By exploiting the operators $W_\Omega$ and $W^*_\Omega$ we can now write the jump formulas
\begin{equation}\label{jump}
w^\pm_\Omega[\psi]_{|\partial\Omega}=\pm\frac{1}{2}\psi+W_\Omega[\psi]\quad\text{ and }\quad\nu_\Omega\cdot\nabla v^\pm_\Omega[\psi]_{|\partial\Omega}=\mp\frac{1}{2}\phi+W^*_\Omega[\psi]
\end{equation}
which hold for all continuous function $\psi\in C^0(\partial\Omega)$ (cf., {\it e.g.},  Folland \cite[Chap.~3]{Fo95}). In addition, if $\psi\in C^{1,\alpha}(\partial\Omega)$, then we have
\begin{equation}\label{nojump}
\nu_\Omega\cdot\nabla w^+_\Omega[\psi]_{|\partial\Omega}=\nu_\Omega\cdot\nabla w^-_\Omega[\psi]_{|\partial\Omega}\,.
\end{equation}

In the following Lemma \ref{ker} we describe the null-spaces  $\mathrm{Ker}(\pm\frac{1}{2}I_\Omega+W^*_\Omega)$  and $\mathrm{Ker}(\pm\frac{1}{2}I_\Omega+W_\Omega)$ of the operators $\pm\frac{1}{2}I_\Omega+W^*_\Omega$ and $\pm\frac{1}{2}I_\Omega+W_\Omega$ in $L^2(\partial\Omega)$.
To do so we exploit the following notation: if $\mathcal{X}$ is a subspace of $L^1(\partial\Omega)$ then we denote by $\mathcal{X}_0$  the subspace of $\mathcal{X}$ consisting of the functions which have $0$ integral mean. For a proof of Lemma \ref{ker} we refer, {\it e.g.},  to Folland \cite[Chap.~3]{Fo95}.

\begin{lem}\label{ker} Let $\Omega_1$, \dots, $\Omega_N$ be the bounded connected components of $\Omega$ and $\Omega^-_0$, $\Omega^-_1$, \dots, $\Omega^-_M$  be the connected components of $\mathbb{R}^n\setminus\mathrm{cl}\Omega$. Assume that $\Omega^-_1$, \dots, $\Omega^-_M$ are bounded and $\Omega^-_0$ is unbounded. Then the following statements hold.
\begin{enumerate}
\item[(i)] The map from $\mathrm{Ker}(\frac{1}{2}I_\Omega+W^*_\Omega)$ to $\mathrm{Ker}(\frac{1}{2}I_\Omega+W_\Omega)$ which takes $\mu$ to $v_\Omega[\mu]_{|\partial\Omega}$ is bijective.
\item[(ii)] The map from $\mathrm{Ker}(-\frac{1}{2}I_\Omega+W^*_\Omega)_0$ to $\mathrm{Ker}(-\frac{1}{2}I_\Omega+W_\Omega)$ which takes $\mu$ to $v_\Omega[\mu]_{|\partial\Omega}$ is one-to-one. If $n\ge 3$, then the map from $\mathrm{Ker}(-\frac{1}{2}I_\Omega+W^*_\Omega)$ to $\mathrm{Ker}(-\frac{1}{2}I_\Omega+W_\Omega)$ which takes $\mu$ to $v_\Omega[\mu]_{|\partial\Omega}$ is bijective.
\item[(iii)]    $\mathrm{Ker}(\frac{1}{2}I_\Omega+W_\Omega)$ consists of the functions from $\partial\Omega$ to $\mathbb{R}$ which are constant on $\partial\Omega^-_j$ for all $j\in\{1,\dots,M\}$ and which are identically equal to $0$ on $\partial\Omega^-_0$.
\item[(iv)] $\mathrm{Ker}(-\frac{1}{2}I_\Omega+W_\Omega)$ consists of the functions from $\partial\Omega$ to $\mathbb{R}$ which are constant on $\partial\Omega_j$ for all $j\in\{1,\dots,N\}$.
\item[(v)]  If $\phi\in\mathrm{Ker}(\frac{1}{2}I_\Omega+W^*_{\Omega})$ and $\int_{\partial\Omega}\phi\psi\, d\sigma=0$ for all  $\psi\in\mathrm{Ker}(\frac{1}{2}I_\Omega+W_{\Omega})$, then $\phi=0$.
\item[(vi)]  If $\phi\in\mathrm{Ker}(-\frac{1}{2}I_\Omega+W^*_{\Omega})$  and $\int_{\partial\Omega}\phi\psi\, d\sigma=0$ for all  $\psi\in\mathrm{Ker}(-\frac{1}{2}I_\Omega+W_{\Omega})$, then $\phi=0$.
\end{enumerate}
\end{lem}

Finally, we have the following technical Lemma  \ref{ilw}.

\begin{lem}\label{ilw}
 Let $\tau\in]-1,1[$.  Then  $\frac{1}{2}I_{\Omega}+\tau W^*_{\Omega}$   is an isomorphism from $L^2(\partial\Omega)$ to itself, from $C^0(\partial\Omega)$ to itself, and from $C^{0,\alpha}(\partial\Omega)$ to itself.  
\end{lem}
\proof We first prove that $\frac{1}{2}I_{\Omega}+\tau W^*_{\Omega}$   is an isomorphism from $L^2(\partial\Omega)$ to itself. To do so we observe that  $\tau W^*_{\Omega}$ is compact from from $L^2(\partial\Omega)$ to itself and thus $\frac{1}{2}I_{\Omega}+\tau W^*_{\Omega}$   is a Fredholm operator of index $0$ from $L^2(\partial\Omega)$ to itself. Accordingly, it suffices to show that $\frac{1}{2}I_{\Omega}+\tau W^*_{\Omega}$ is one-to-one.  A fact which can be verified by arguing as in Escauriaza {\it et al.}~\cite[\S3 (i)]{EsFaVe92}. To prove that $\frac{1}{2}I_{\Omega}+\tau W^*_{\Omega}$  is invertible from $C^0(\partial\Omega)$ to itself we observe that $\frac{1}{2}I_{\Omega}+\tau W^*_{\Omega}$ is continuous from $C^0(\partial\Omega)$ to itself (because $W^*_{\Omega}$ has a weak singularity). Moreover, if $\eta\in L^2(\partial\Omega)$ and $(\frac{1}{2}I_{\Omega}+\tau W^*_{\Omega})\eta\in C^0(\partial\Omega)$ then  Lemma \ref{regularity} ensures that $\eta\in C^{0}(\partial\Omega)$. Similarly, to prove that $\frac{1}{2}I_{\Omega}+\gamma W^*_{\Omega}$ is invertible from $C^{0,\alpha}(\partial\Omega)$ to itself we observe that $\frac{1}{2}I_{\Omega}+\tau W^*_{\Omega}$ is continuous from $C^{0,\alpha}(\partial\Omega)$ to itself and that $(\frac{1}{2}I_{\Omega}+\tau W^*_{\Omega})\eta\in C^{0,\alpha}(\partial\Omega)$ implies $\eta\in C^{0,\alpha}(\partial\Omega)$ for all $\eta\in L^2(\partial\Omega)$ by  Lemma \ref{regularity}.\qed

%Finally, if $\Omega$ is an open bounded subset of $\mathbb{R}^2$ of class $C^{1,\alpha}$ and $\mathcal{X}$ is a subspace of  $L^1(\partial\Omega)$,  then we denote by   $\mathcal{X}_0$ the subspace of $\mathcal{X}$  consisting of those functions $f$ such that $\int_{\partial\Omega} f\,d\sigma=0$.

\section{Existence results for problem \eqref{nlprob}}\label{sec.existence}

We prove in this section our main Theorems \ref{existence}.

As a first step we deduce in the following Lemma \ref{representation} a representation for a pair of harmonic functions in $C^{1,\alpha}(\mathrm{cl}\Omega^o\setminus\Omega^i)\times C^{1,\alpha}(\mathrm{cl}\Omega^i)$ in terms of a suitable combination of layer potential. We will exploit the following notation: if $\Omega$ is an open bounded subset of $\mathbb{R}^n$, $k\in\mathbb{N}$, and $\beta\in[0,1[$, then we denote by $C^{k,\beta}_\mathrm{harm}(\mathrm{cl}\Omega)$ the subspace of $C^{k,\beta}(\mathrm{cl}\Omega)$ defined by
\begin{equation}\label{Charm}
C^{k,\beta}_\mathrm{harm}(\mathrm{cl}\Omega)\equiv\left\{\phi\in C^{k,\beta}(\mathrm{cl}\Omega)\,:\; \Delta\phi=0\text{ in }\Omega\right\}\,.
\end{equation}

\begin{lem}\label{representation}
The map from $C^{1,\alpha}(\partial\Omega^o)\times C^{1,\alpha}(\partial\Omega^i)\times C^{0,\alpha}(\partial\Omega^i)$ to $C^{1,\alpha}_\mathrm{harm}(\mathrm{cl}\Omega^o\setminus\Omega^i)\times C^{1,\alpha}_\mathrm{harm}(\mathrm{cl}\Omega^i)$ which takes $(\mu^o,\mu,\eta)$ to the pair $(u^o[\mu^o,\mu,\eta], u^i[\mu^o,\mu,\eta])$ given by 
\[
u^o[\mu^o,\mu,\eta]\equiv(w^+_{\Omega^o}[\mu^o]+w^-_{\Omega^i}[\mu]+v^-_{\Omega^i}[\eta])_{|\mathrm{cl}\Omega^o\setminus\Omega^i}\,,\quad u^i[\mu^o,\mu,\eta]\equiv w^+_{\Omega^i}[\mu]
\] is bijective.
\end{lem}
\proof 
The map is well defined. Indeed $(u^o[\mu^o,\mu,\eta], u^i[\mu^o,\mu,\eta])\in C^{1,\alpha}(\mathrm{cl}\Omega^o\setminus\Omega^i)\times C^{1,\alpha}(\mathrm{cl}\Omega^i)$ and $\Delta u^o[\mu^o,\mu,\eta]=0$, $\Delta u^i[\mu^o,\mu,\eta]=0$ for all $(\mu^o,\mu,\eta)\in C^{1,\alpha}(\partial\Omega^o)\times C^{1,\alpha}(\partial\Omega^i)\times C^{0,\alpha}(\partial\Omega^i)$ (cf.~Section \ref{sec.preliminaries}). We now show that it is bijective. We take a pair of harmonic functions $(\phi^o,\phi^i)$ in $C^{1,\alpha}_\mathrm{harm}(\mathrm{cl}\Omega^o\setminus\Omega^i)\times C^{1,\alpha}_\mathrm{harm}(\mathrm{cl}\Omega^i)$ and we prove that there exists unique $(\mu^o,\mu,\eta)\in C^{1,\alpha}(\partial\Omega^o)\times C^{1,\alpha}(\partial\Omega^i)\times C^{0,\alpha}(\partial\Omega^i)$ such that $(u^o[\mu^o,\mu,\eta], u^i[\mu^o,\mu,\eta])=(\phi^o,\phi^i)$. By the standard properties of the double layer potential there exists a unique $\mu\in C^{1,\alpha}(\partial\Omega^i)$ such that $w^+_{\Omega^i}[\mu]=\phi^i$ (cf.~\eqref{jump} and Lemma \ref{ker} (iii)). Then we have to show that there exists unique $(\mu^o,\eta)\in C^{1,\alpha}(\partial\Omega^o)\times C^{0,\alpha}(\partial\Omega^i)$ such that 
\begin{equation}\label{representation.eq1}
(w^+_{\Omega^o}[\mu^o]+v^-_{\Omega^i}[\eta])_{|\mathrm{cl}\Omega^o\setminus\Omega^i}=\phi^o-w^-_{\Omega^i}[\mu]_{|\mathrm{cl}\Omega^o\setminus\Omega^i}\,.
\end{equation} 
Let $\psi^o\equiv \phi^o_{|\partial\Omega^o}-w^-_{\Omega^i}[\mu]_{|\partial\Omega^o}$ and $\psi^i\equiv \nu_{\Omega^i}\cdot\nabla \phi^o_{|\partial\Omega^i}-\nu_{\Omega^i}\cdot\nabla w^-_{\Omega^i}[\mu]_{|\partial\Omega^i}$. Then $\psi^o\in C^{1,\alpha}(\partial\Omega^o)$, $\psi^i\in C^{0,\alpha}(\partial\Omega^o)$, and equation  \eqref{representation.eq1} is equivalent to 
\begin{equation}\label{representation.eq2}
\begin{split}
&(\frac{1}{2}I_{\Omega^o}+W_{\Omega^o})\mu^o+v^-_{\Omega^i}[\eta]_{|\partial\Omega^o}=\psi^o\,,\\
&(\frac{1}{2}I_{\Omega^i}+W^*_{\Omega^i})\eta+\nu_{\Omega^i}\cdot\nabla w^+_{\Omega^o}[\mu^o]_{|\partial\Omega^i}=\psi^i
\end{split}
\end{equation}
 by the uniqueness of the classical solution of the Neumann-Dirichlet mixed boundary value problem  (see also \eqref{jump}).
 By Lemmas \ref{fredC} and \ref{ker} the operator which takes $(\mu^o,\eta)$ to 
 $((\frac{1}{2}I_{\Omega^o}+W_{\Omega^o})\mu^o, (\frac{1}{2}I_{\Omega^i}+W^*_{\Omega^i})\eta)$ is a linear isomorphism from $C^{1,\alpha}(\partial\Omega^o)\times C^{0,\alpha}(\partial\Omega^i)$ to itself. Moreover, by the properties of the integral operators with real analytic kernels and no singularities, the operator which takes $(\mu^o,\eta)$ to $(v^-_{\Omega^i}[\eta]_{|\partial\Omega^o},\nu_{\Omega^i}\cdot\nabla w^+_{\Omega^o}[\mu^o]_{|\partial\Omega^i})$ is compact from $C^{1,\alpha}(\partial\Omega^o)\times C^{0,\alpha}(\partial\Omega^i)$ to itself. Hence,  the operator which takes $(\mu^o,\eta)$ to $((\frac{1}{2}I_{\Omega^o}+W_{\Omega^o})\mu^o+v^-_{\Omega^i}[\eta]_{|\partial\Omega^o}, (\frac{1}{2}I_{\Omega^i}+W^*_{\Omega^i})\eta+\nu_{\Omega^i}\cdot\nabla w^+_{\Omega^o}[\mu^o]_{|\partial\Omega^i})$ is a compact perturbation of an isomorphism and therefore it is a Fredholm operator of index $0$ from $C^{1,\alpha}(\partial\Omega^o)\times C^{0,\alpha}(\partial\Omega^i)$ to itself. Thus, to complete the proof it suffices to show that equation \eqref{representation.eq2} with $(\psi^o,\psi^i)=(0,0)$ implies $(\mu^o,\eta)=(0,0)$. If $((\frac{1}{2}I_{\Omega^o}+W_{\Omega^o})\mu^o+v^-_{\Omega^i}[\eta]_{|\partial\Omega^o}, (\frac{1}{2}I_{\Omega^i}+W^*_{\Omega^i})\eta+\nu_{\Omega^i}\cdot\nabla w^+_{\Omega^o}[\mu^o]_{|\partial\Omega^i})=(0,0)$,  then by the jump properties \eqref{jump} and  by the uniqueness of the classical solution of the Neumann-Dirichlet mixed problem one deduces that $(w^+_{\Omega^o}[\mu^o]+v^-_{\Omega^i}[\eta])_{|\mathrm{cl}\Omega^o\setminus\Omega^i}=0$. Hence $w^+_{\Omega^o}[\mu^o]+v^+_{\Omega^i}[\eta]=0$ in $\mathrm{cl}\Omega^i$ by the uniqueness of the classical solution of the  Dirichlet  problem in $\Omega^i$ and by the continuity of $(w^+_{\Omega^o}[\mu^o]+v_{\Omega^i}[\eta])_{|\mathrm{cl}\Omega^o}$ (cf.~Section \ref{sec.preliminaries}). Then by \eqref{jump}  we have
 \[
 \begin{split}
\eta&=\nu_{\Omega^i}\cdot\nabla v^-_{\Omega^i}[\eta]_{|\partial\Omega^i}-\nu_{\Omega^i}\cdot\nabla v^+_{\Omega^i}[\eta]_{|\partial\Omega^i}\\
&=\nu_{\Omega^i}\cdot\nabla(w^+_{\Omega^o}[\mu^o]+v^-_{\Omega^i}[\eta])_{|\partial\Omega^i}-\nu_{\Omega^i}\cdot\nabla(w^+_{\Omega^o}[\mu^o]+v^+_{\Omega^i}[\eta])_{|\partial\Omega^i}=0\,.
 \end{split}
 \]
  By \eqref{representation.eq2} it follows that  $(\frac{1}{2}I_{\Omega^o}+W_{\Omega^o})\mu^o=0$ and thus  $\mu^o=0$ by Lemma \ref{ker} (iii). Our proof is now complete.
\qed

\vspace{\baselineskip}

In the following Lemma \ref{J} we introduce an auxiliary operator which we denote by $J$. In the sequel we will denote the inverse of an 
invertible map $f$ with $f^{(-1)}$, as opposed to the 
reciprocal of a function $g$ which will be denoted with $g^{-1}$.

\begin{lem}\label{J}
We define
\[
J[\eta]\equiv(\frac{1}{2}I_{\Omega^i}+W^*_{\Omega^i})\eta-\nu_{\Omega^i}\cdot\nabla w^+_{\Omega^o}\Bigl[(\frac{1}{2}I_{\Omega^o}+W_{\Omega^o})^{(-1)}v_{\Omega^i}[\eta]_{|\partial\Omega^o}\Bigr]_{|\partial\Omega^i}
\] for all $\eta\in L^2(\partial\Omega^i)$. Then the map which takes $\eta$ to $J[\eta]$ is an isomorphism from $L^2(\partial\Omega^i)$ to itself,  from $C^0(\partial\Omega^i)$ to itself, and from $C^{0,\alpha}(\partial\Omega^i)$ to itself.
\end{lem}
\proof By the properties of integral operators with real analytic kernels and no singularity, by the invertibility of $\frac{1}{2}I_{\Omega^o}+W_{\Omega^o}$ in $C^{1,\alpha}(\partial\Omega^o)$ (cf.~Lemmas \ref{fredC} and \ref{ker}), and by the continuity of the map $w^+_{\Omega^o}[\cdot]$ from $C^{1,\alpha}(\partial\Omega^o)$ to $C^{1,\alpha}(\mathrm{cl}\Omega^o)$ (cf., {\it e.g.}, Miranda \cite{Mi65}), one deduces that the operator which takes $\eta$ to  
\begin{equation}\label{J.eq1}
\nu_{\Omega^i}\cdot\nabla w^+_{\Omega^o}\Bigl[(\frac{1}{2}I_{\Omega^o}+W_{\Omega^o})^{(-1)}v_{\Omega^i}[\eta]_{|\partial\Omega^o}\Bigr]_{|\partial\Omega^i}
\end{equation}
is continuous from $L^2(\partial\Omega^i)$ to $C^{0,\alpha}(\partial\Omega^i)$. Then, by the compactness of $W^*_{\Omega^i}$ in $L^2(\partial\Omega^i)$ it follows that $J$ is Fredholm operator of index $0$ from $L^2(\partial\Omega^i)$ to itself. Thus, to show that $J$ is invertible from $L^2(\partial\Omega^i)$ to itself it suffices to prove that $J[\eta]=0$  implies $\eta=0$. If $\eta\in L^2(\partial\Omega^i)$ and $J[\eta]=0$, then $(\frac{1}{2}I_{\Omega^i}+W^*_{\Omega^i})\eta\in C^{0,\alpha}(\partial\Omega^i)$ by the membership of \eqref{J.eq1} in $C^{0,\alpha}(\partial\Omega^i)$. It follows that $\eta\in C^{0,\alpha}(\partial\Omega^i)$ (cf.~Lemma \ref{regularity}). Then, by setting $\mu^o\equiv -(\frac{1}{2}I_{\Omega^o}+W_{\Omega^o})^{(-1)}v_{\Omega^i}[\eta]_{|\partial\Omega^o}$ and by exploiting equality \eqref{jump} we verify that $u^o[\mu^o,0,\eta]_{|\partial\Omega^o}=0$ and $\nu_{\Omega^i}\cdot\nabla u^o[\mu^o,0,\eta]_{|\partial\Omega^o}=0$, where $u^o[\mu^o,0,\eta]$ is defined as in Lemma \ref{representation}. Accordingly  $u^o[\mu^o,0,\eta]=0$ by the uniqueness of the solution of the mixed boundary value problem. Since $u^i[\mu^o,0,\eta]=w^+_{\Omega^i}[0]=0$, Lemma \ref{representation} implies that $\eta=0$. 

To prove that $J$ is invertible from $C^0(\partial\Omega^i)$ to itself we observe that $J$ is continuous from $C^0(\partial\Omega^i)$ to itself (because $W^*_{\Omega^i}$ has a weak singularity). Moreover, if $\eta\in L^2(\partial\Omega^i)$ and $J[\eta]\in C^0(\partial\Omega^i)$ then $(\frac{1}{2}I_{\Omega^i}+W^*_{\Omega^i})\eta\in C^{0}(\partial\Omega^i)$ by the membership of \eqref{J.eq1} in $C^{0,\alpha}(\partial\Omega^i)$. Thus Lemma \ref{regularity} ensures that $\eta\in C^{0}(\partial\Omega^i)$. 

Similarly, to prove that $J$ is invertible from $C^{0,\alpha}(\partial\Omega^i)$ to itself we observe that $J$ is continuous from $C^{0,\alpha}(\partial\Omega^i)$ to itself and that $J[\eta]\in C^{0,\alpha}(\partial\Omega^i)$ implies $\eta\in C^{0,\alpha}(\partial\Omega^i)$ for all $\eta\in L^2(\partial\Omega^i)$.\qed

\vspace{\baselineskip}

Then we have the following Lemma \ref{solution} where we rewrite problem \eqref{nlprob} into an equivalent system of boundary integral equations.

%To shorten our notation we also  introduce the  operator $M^o$ from $(L^2(\partial\Omega^i))^3$ to $L^2(\partial\Omega^i)$ defined by
%\begin{equation}\label{M}
%M^o[\psi^o,\psi,\phi]\equiv (\frac{1}{2}I_{\Omega^o}+W_{\Omega^o})^{(-1)}(\psi^o-w^-_{\Omega^i}[\psi]_{|\partial\Omega^o}-v_{\Omega^i}[\phi]_{|\partial\Omega^o})\qquad\forall(\psi^o,\psi,\phi)\in (L^2(\partial\Omega^i))^3\,.
%\end{equation}

\begin{lem}\label{solution} Let condition \eqref{assumption0} hold.
Let $(\mu^o,\mu,\eta)\in C^{1,\alpha}(\partial\Omega^o)\times C^{1,\alpha}(\partial\Omega^i)\times C^{0,\alpha}(\partial\Omega^i)$.  Then $(u^o[\mu^o,\mu,\eta], u^i[\mu^o,\mu,\eta])$ is a solution of \eqref{nlprob} if and only if 
\begin{equation}\label{integral}
\begin{split}
\mu^o&=(\frac{1}{2}I_{\Omega^o}+W_{\Omega^o})^{(-1)}(f^o-w^-_{\Omega^i}[\mu]_{|\partial\Omega^o}-v_{\Omega^i}[\eta]_{|\partial\Omega^o})\,,\\
\mu&=(\frac{1}{2}I_{\Omega^i}+W_{\Omega^i})^{(-1)}\left[ (I_{\Omega^i}+\mathcal{F}_F)^{(-1)}\left(w^+_{\Omega^o}[\mu^o]_{|\partial\Omega^i}+v_{\Omega^i}[\eta]_{|\partial\Omega^i}+2W_{\Omega^i}\mu\right)\right]\,,\\
\eta&=J^{(-1)}\bigg[\mathcal{F}_G\circ(I_{\Omega^i}+\mathcal{F}_F)^{(-1)}\left(w^+_{\Omega^o}[\mu^o]_{|\partial\Omega^i}+v_{\Omega^i}[\eta]_{|\partial\Omega^i}+2W_{\Omega^i}\mu\right)\\
&\qquad\qquad -\nu_{\Omega^i}\cdot\nabla w^+_{\Omega^o}[(\frac{1}{2}I_{\Omega^o}+W_{\Omega^o})^{(-1)}(f^o-w^-_{\Omega^i}[\mu]_{|\partial\Omega^o})]_{|\partial\Omega^i}\bigg]\,.\\
\end{split}
\end{equation}
\end{lem}
\proof Note that $\nu_{\Omega^i}\cdot\nabla w^-_{\Omega^i}[\mu](x)-\nu_{\Omega^i}\cdot\nabla w^+_{\Omega^i}[\mu](x)=0$ by the membership of $\mu$ in $C^{1,\alpha}(\partial\Omega^i)$ (cf.~\eqref{nojump}). Then the validity of the statement is a consequence of Lemma \ref{representation}, of the jump properties of single and double layer potentials (cf.~\eqref{jump}), of the invertibility of $(\frac{1}{2}I_{\Omega^o}+W_{\Omega^o})$ in $C^{1,\alpha}(\partial\Omega^o)$,  of the invertibility of  $(\frac{1}{2}I_{\Omega^i}+W_{\Omega^i})$ and  $J$ in $L^2(\partial\Omega^i)$ (cf.~Lemmas \ref{ker} and \ref{J}), and of condition \eqref{assumption0}.\qed 

\vspace{\baselineskip}

In Proposition \ref{existenceC} below we prove the existence of a solution $(\tilde\mu^o,\tilde\mu,\tilde\eta)$ in $C^{1,\alpha}(\partial\Omega^i)\times C^0(\partial\Omega^i)\times C^0(\partial\Omega^i)$ of the system of equations in \eqref{integral}. To do so we exploit the Leray-Schauder principle  which is stated in the following Theorem \ref{topological} and which follows by the invariance of the Leray-Schauder topological degree (for a proof see, {\it e.g.}, Gilbarg and Trudinger \cite[Theorem 11.3]{GiTr01}).

\begin{thm}[Leray-Schauder principle]\label{topological} Let $\mathcal{X}$ be a Banach space. Let $T$ be a continuous (nonlinear) operator from $\mathcal{X}$ to itself which maps bounded sets to sets with a compact closure. If there exists a constant $M\in]0,+\infty[$ such that $\|x\|_\mathcal{X}\le M$ for all $(x,t)\in\mathcal{X}\times[0,1]$ satisfying $x = tT(x)$, then $T$ has at least one fixed point $x\in\mathcal{X}$ such that $\|x\|_\mathcal{X}\le M$.
\end{thm}

In order to apply this principle, we   introduce in the following Lemma an elementary consequence of conditions \eqref{assumption0} and \eqref{assumption1}.

\begin{lem}\label{C}
Let conditions \eqref{assumption0} and \eqref{assumption1} hold. Then there exist $C_1,C_2, C_3, C_4\in]0,+\infty[$ such that 
\begin{equation}\label{C.eq1}
\|(I_{\Omega^i}+\mathcal{F}_F)^{(-1)}f\|_{C^0(\partial\Omega^i)}\le C_1(C_2+\|f\|_{C^0(\partial\Omega^i)})^{1/\delta_1}
\end{equation}
and
\begin{equation}\label{C.eq2}
\|\mathcal{F}_G\circ (I_{\Omega^i}+\mathcal{F}_F)^{(-1)}f\|_{C^0(\partial\Omega^i)}\le C_3(C_4+\|f\|_{C^0(\partial\Omega^i)})^{\delta_2}
\end{equation}
for all functions $f\in C^0(\partial\Omega^i)$.
\end{lem}
\proof 
To prove \eqref{C.eq1} we observe that the first inequality in \eqref{assumption1} implies that there exist $c_1^*,c_2^*\in]0,+\infty[$ such that $|t+F(x,t)|\ge c_1^*|t|^{\delta_1}-c_2^*$ for all $(x,t)\in\partial\Omega^i\times\mathbb{R}$. Thus we have  $\|(I_{\Omega^i}+\mathcal{F}_F)g\|_{C^0(\partial\Omega^i)}\ge c_1^*\|g\|_{C^0(\partial\Omega^i)}^{\delta_1}-c_2^*$ for all $g\in C^0(\partial\Omega^i)$ and the validity of \eqref{C.eq1} follows by  taking $g= (I_{\Omega^i}+\mathcal{F}_F)^{(-1)}f $. To prove \eqref{C.eq2} we observe that  the second inequality in \eqref{assumption1} implies that there exist $c_3^*,c_4^*\in]0,+\infty[$ such that $|G(x,t)|\le c_3^*(c_4^*+|t+F(x,t)|)^{\delta_2}$ for all $(x,t)\in\partial\Omega^i\times\mathbb{R}$. Then we have $\|\mathcal{F}_Gg\|_{C^0(\partial\Omega^i)}\le c_3^*(c_4^*+\|(I_{\Omega^i}+\mathcal{F}_F)g\|_{C^0(\partial\Omega^i)})^{\delta_2}$ for all $g\in C^0(\partial\Omega^i)$ and the validity of \eqref{C.eq2} follows by condition \eqref{assumption0} and by  taking $g= (I_{\Omega^i}+\mathcal{F}_F)^{(-1)}f $. 
\qed

\vspace{\baselineskip}

Then we have the following.

\begin{prop}\label{existenceC} Let conditions \eqref{assumption0} and \eqref{assumption1} hold. Then the nonlinear system \eqref{integral} has at least one solution $(\tilde\mu^o,\tilde\mu,\tilde\eta)$ in $C^{1,\alpha}(\partial\Omega^o)\times C^0(\partial\Omega^i)\times C^0(\partial\Omega^i)$. 
\end{prop}
\proof We plan to apply Theorem \ref{topological} with $\mathcal{X}=C^{1,\alpha}(\partial\Omega^o)\times C^0(\partial\Omega^i)\times C^0(\partial\Omega^i)$ and $T\equiv(T^o,T_1,T_2)$ given by
\begin{equation}\label{T}
\begin{split}
&T^o(\tilde\mu^o,\tilde\mu,\tilde\eta)\equiv (\frac{1}{2}I_{\Omega^o}+W_{\Omega^o})^{(-1)}(f^o-w^-_{\Omega^i}[\tilde\mu]_{|\partial\Omega^o}-v_{\Omega^i}[\tilde\eta]_{|\partial\Omega^o})\\
&T_1(\tilde\mu^o,\tilde\mu,\tilde\eta)\\
&\quad \equiv(\frac{1}{2}I_{\Omega^i}+W_{\Omega^i})^{(-1)}\left[ (I_{\Omega^i}+\mathcal{F}_F)^{(-1)}\left(w^+_{\Omega^o}[\tilde\mu^o]_{|\partial\Omega^i}+v_{\Omega^i}[\tilde\eta]_{|\partial\Omega^i}+2W_{\Omega^i}\tilde\mu\right)\right],\\
&T_2(\tilde\mu^o,\tilde\mu,\tilde\eta)\equiv J^{(-1)}\bigg[\mathcal{F}_G\circ(I_{\Omega^i}+\mathcal{F}_F)^{(-1)}\left(w^+_{\Omega^o}[\tilde\mu^o]_{|\partial\Omega^i}+v_{\Omega^i}[\tilde\eta]_{|\partial\Omega^i}+2W_{\Omega^i}\tilde\mu\right)\\
&\qquad\qquad\qquad -\nu_{\Omega^i}\cdot\nabla w^+_{\Omega^o}[(\frac{1}{2}I_{\Omega^o}+W_{\Omega^o})^{(-1)}(f^o-w^-_{\Omega^i}[\tilde\mu]_{|\partial\Omega^o})]_{|\partial\Omega^i}\bigg]\,,\\
\end{split}
\end{equation}  for all $(\tilde\mu^o,\tilde\mu,\tilde\eta)\in C^{1,\alpha}(\partial\Omega^o)\times C^0(\partial\Omega^i)\times C^0(\partial\Omega^i)$. We first verify that $T$ is continuous from $C^{1,\alpha}(\partial\Omega^o)\times C^0(\partial\Omega^i)\times C^0(\partial\Omega^i)$ to itself and maps bounded sets to sets with compact closure. To do so, we consider separately $T^o$, $T_1$ and $T_2$. By Lemmas  \ref{fredC} and \ref{ker} one deduces that $\frac{1}{2}I_{\Omega^o}+W_{\Omega^o}$  is an isomorphism from $C^{1,\alpha}(\partial\Omega^o)$ to itself. In particular, $(\frac{1}{2}I_{\Omega^o}+W_{\Omega^o})^{(-1)}$ is continuous from  $C^{1,\alpha}(\partial\Omega^o)$ to itself. Moreover, by the properties of integral operators with real analytic kernel and no singularities $w^-_{\Omega^i}[\cdot]_{|\partial\Omega^o}$ and $v_{\Omega^i}[\cdot]_{|\partial\Omega^o}$ are compact from $C^0(\partial\Omega^i)$ to $C^{1,\alpha}(\partial\Omega^o)$. It follows that $T^o$ is continuous from $C^{1,\alpha}(\partial\Omega^o)\times C^0(\partial\Omega^i)\times C^0(\partial\Omega^i)$ to $C^{1,\alpha}(\partial\Omega^o)$ and maps bounded sets to sets with compact closure. We now consider $T_1$. By  Lemmas  \ref{fredC} and \ref{ker} one verifies that $\frac{1}{2}I_{\Omega^i}+W_{\Omega^i}$ is an isomorphism from $C^0(\partial\Omega^i)$ to itself and thus
$(\frac{1}{2}I_{\Omega^i}+W_{\Omega^i})^{(-1)}$ is continuous from $C^0(\partial\Omega^i)$ to itself. By assumption  \eqref{assumption0}  the map $(I_{\Omega^i}+\mathcal{F}_F)^{(-1)}$ is  continuous from $C^0(\partial\Omega^i)$ to itself. Then, by the properties of integral operators with real analytic kernel and no singularities $w^+_{\Omega^o}[\cdot]_{|\partial\Omega^i}$ is compact from $C^{1,\alpha}(\partial\Omega^o)$ to $C^0(\partial\Omega^i)$. By the mapping properties of the single layer potential  (cf., {\it e.g.}, Kress \cite[Thm.~2.22]{Kr89}, see also Miranda \cite[Chap.~II, \S14, III]{Mi70}), $v_{\Omega^i}[\cdot]_{|\partial\Omega^i}$ is  compact from $C^0(\partial\Omega^i)$ to itself. By Lemma \ref{schauder}, $W_{\Omega^i}$ is compact from $C^0(\partial\Omega^i)$ to itself. If follows that  
$T_1$ is continuous from $C^{1,\alpha}(\partial\Omega^o)\times C^0(\partial\Omega^i)\times C^0(\partial\Omega^i)$ to $C^{0}(\partial\Omega^i)$ and maps bounded sets to sets with compact closure. Finally we consider $T_2$. By Lemma \ref{J} the operator $J^{(-1)}$ is continuous from $C^0(\partial\Omega^i)$ to itself. By the continuity of $G$ and by condition \eqref{assumption0}, the map  $\mathcal{F}_G\circ(I_{\Omega^i}+\mathcal{F}_F)^{(-1)}$ is continuous from $C^0(\partial\Omega^i)$ to itself.  By the mapping properties of the single layer potential  (cf., {\it e.g.}, Kress \cite[Thm.~2.22]{Kr89}, see also Miranda \cite[Chap.~II, \S14, III]{Mi70}), $v_{\Omega^i}[\cdot]_{|\partial\Omega^i}$ is  compact from $C^0(\partial\Omega^i)$ to itself.  By Lemma \ref{schauder},   $W_{\Omega^i}$ is compact from $C^0(\partial\Omega^i)$ to itself.  By the properties of integral operators with real analytic kernel and no singularities and by the continuity of $(\frac{1}{2}I_{\Omega^o}+W_{\Omega^o})^{(-1)}$ from $C^{1,\alpha}(\partial\Omega^o)$ to itself,  the map $w^+_{\Omega^o}[\cdot]_{|\partial\Omega^i}$ is compact from $C^{1,\alpha}(\partial\Omega^o)$ to $C^0(\partial\Omega^i)$ and the map $\nu_{\Omega^i}\cdot\nabla w^+_{\Omega^o}[(\frac{1}{2}I_{\Omega^o}+W_{\Omega^o})^{(-1)}w^-_{\Omega^i}[\cdot]_{|\partial\Omega^o}]_{|\partial\Omega^i}$ is compact from $C^0(\partial\Omega^i)$ to itself. Accordingly $T_2$ is continuous from $C^{1,\alpha}(\partial\Omega^o)\times C^0(\partial\Omega^i)\times C^0(\partial\Omega^i)$ to $C^{0}(\partial\Omega^i)$ and maps bounded sets to sets with compact closure.

 Now let $t\in[0,1]$ and assume that $(\tilde\mu^o,\tilde\mu,\tilde\eta)=tT(\tilde\mu^o,\tilde\mu,\tilde\eta)$. We show that there exists a constant $M\in]0,+\infty[$ (which does not depend on $t$) such that 
 \begin{equation}\label{existenceC.eq-1}
\|\tilde\mu^o\|_{C^{1,\alpha}(\partial\Omega^o)}+\|\tilde\mu\|_{C^{0}(\partial\Omega^i)}+\|\tilde\eta\|_{C^{0}(\partial\Omega^i)}\le M\,.
 \end{equation} 
By equality $(\tilde\mu^o,\tilde\mu,\tilde\eta)=tT(\tilde\mu^o,\tilde\mu,\tilde\eta)$ we have that
 \begin{equation}
 \label{existenceC.eq0}
 \begin{split}
&\|\tilde\mu^o\|_{C^{1,\alpha}(\partial\Omega^o)}\le \|T^o(\tilde\mu^o,\tilde\mu,\tilde\eta)\|_{C^{1,\alpha}(\partial\Omega^o)}\,,\\
& \|\tilde\mu\|_{C^0(\partial\Omega^i)}\le \|T_1(\tilde\mu^o,\tilde\mu,\tilde\eta)\|_{C^0(\partial\Omega^i)}\,,\\
&\|\tilde\eta\|_{C^0(\partial\Omega^i)}\le \|T_2(\tilde\mu^o,\tilde\mu,\tilde\eta)\|_{C^0(\partial\Omega^i)}\,.
\end{split}
\end{equation}
By the first inequality of \eqref{existenceC.eq0} we deduce that there exists a constant $m_1\in]0,+\infty[$ which depends only on the norm of the bounded linear operator $(\frac{1}{2}I_{\Omega^o}+W_{\Omega^o})^{(-1)}$ from $C^{1,\alpha}(\partial\Omega^o)$ to itself, on $\|f^o\|_{C^{1,\alpha}(\partial\Omega^o)}$, and on the norm of the linear bounded operators $w^-_{\Omega^i}[\cdot]_{|\partial\Omega^o}$ and $v_{\Omega^i}[\cdot]_{|\partial\Omega^o}$ from $C^{0}(\partial\Omega^i)$ to $C^{1,\alpha}(\partial\Omega^o)$, such that
\begin{equation}\label{existenceC.eq1}
\|\tilde\mu^o\|_{C^{1,\alpha}(\partial\Omega^o)}\le m_1(1+\|\tilde\mu\|_{C^{0}(\partial\Omega^i)}+\|\tilde\eta\|_{C^{0}(\partial\Omega^i)})\,.
\end{equation}
By the second inequality of \eqref{existenceC.eq0} we deduce that there exist real  constants $m_2, m_3\in]0,+\infty[$ which depend on the norm of the linear bounded operator $(\frac{1}{2}I_{\Omega^i}+W_{\Omega^i})^{(-1)}$ from $C^0(\partial\Omega^i)$ to itself, on the constants $C_1$ and $C_2$ of Lemma \ref{C}, on the norm of the linear bounded operator $w^+_{\Omega^o}[\cdot]_{|\partial\Omega^i}$ from $C^{1,\alpha}(\partial\Omega^o)$ to $C^0(\partial\Omega^i)$, and on the norm of the linear bounded operators $v_{\Omega^i}[\cdot]_{|\partial\Omega^i}$ and $W_{\Omega^i}$ from $C^0(\partial\Omega^i)$ to itself such that
\begin{equation}\label{existenceC.eq2}
\|\tilde\mu\|_{C^{0}(\partial\Omega^o)}\le m_2(m_3+\|\tilde\mu^o\|_{C^{1,\alpha}(\partial\Omega^o)}+\|\tilde\mu\|_{C^{0}(\partial\Omega^i)}+\|\tilde\eta\|_{C^{0}(\partial\Omega^i)})^{1/\delta_1}\,.
\end{equation}
By the third inequality of \eqref{existenceC.eq0} we deduce that there exist real  constants $m_4, m_5\in]0,+\infty[$ which depend on the norm of the linear bounded operator $J^{(-1)}$ from $C^0(\partial\Omega^i)$ to itself, on the constants $C_3$ and $C_4$ of Lemma \ref{C}, on the norm of the linear bounded operator $w^+_{\Omega^o}[\cdot]_{|\partial\Omega^i}$ from $C^{1,\alpha}(\partial\Omega^o)$ to $C^0(\partial\Omega^i)$, on the norm of the linear bounded operators $v_{\Omega^i}[\cdot]_{|\partial\Omega^i}$ and $W_{\Omega^i}$ from $C^0(\partial\Omega^i)$ to itself,   on the norm of $\nu_{\Omega^i}\cdot\nabla w^+_{\Omega^o}[(\frac{1}{2}I_{\Omega^o}+W_{\Omega^o})^{(-1)}f^o]_{|\partial\Omega^i}$ in $C^0(\partial\Omega^i)$,  and on the norm of the bounded linear operator $\nu_{\Omega^i}\cdot\nabla w^+_{\Omega^o}[(\frac{1}{2}I_{\Omega^o}+W_{\Omega^o})^{(-1)}w^-_{\Omega^i}[\cdot]_{|\partial\Omega^o}]_{|\partial\Omega^i}$ from $C^0(\partial\Omega^i)$ to itself, such that 
\begin{equation}\label{existenceC.eq3}
\begin{split}
&\|\tilde\eta\|_{C^{0}(\partial\Omega^o)}\\
&\le m_4\left[(m_5+\|\tilde\mu^o\|_{C^{1,\alpha}(\partial\Omega^o)}+\|\tilde\mu\|_{C^{0}(\partial\Omega^i)}+\|\tilde\eta\|_{C^{0}(\partial\Omega^i)})^{\delta_2}+1+\|\tilde\mu\|_{C^{0}(\partial\Omega^i)}\right].
\end{split}
\end{equation}
Then, by inequalities \eqref{existenceC.eq1}, \eqref{existenceC.eq2}, and \eqref{existenceC.eq3} one deduces that there exists real constants $M_1, M_2, M_3\in]0,+\infty[$, which depend on $m_1,\dots,m_5$, such that 
\[
\begin{split}
&\|\tilde\mu^o\|_{C^{1,\alpha}(\partial\Omega^o)}+\|\tilde\mu\|_{C^{0}(\partial\Omega^i)}+\|\tilde\eta\|_{C^{0}(\partial\Omega^i)}\\
&\qquad\qquad\qquad \le M_1+M_2(M_3+\|\tilde\mu^o\|_{C^{1,\alpha}(\partial\Omega^o)}+\|\tilde\mu\|_{C^{0}(\partial\Omega^i)}+\|\tilde\eta\|_{C^{0}(\partial\Omega^i)})^{\delta_*}
\end{split}
\]
with $\delta_*\equiv\max\{1/\delta_1,\delta_2\}\in]0,1[$. Then a straightforward calculation shows that inequality \eqref{existenceC.eq-1} holds with $M\equiv \max\left\{1, (M_1+M_2(M_3+1)^{\delta_*})^{1/(1-\delta_*)}\right\}$.
Now the validity of the statement follows by Theorem \ref{topological}.
\qed

\vspace{\baselineskip}

With a further regularity request on $F$ and $G$ we can find a solution of \eqref{integral} in $C^{1,\alpha}(\partial\Omega^o)\times C^{0,\alpha}(\partial\Omega^i)\times C^{0,\alpha}(\partial\Omega^i)$. 

\begin{prop}\label{C0a}
Let conditions \eqref{assumption0}, \eqref{assumption1}, and \eqref{assumption2} hold.  Then the nonlinear system \eqref{integral} has at least one solution $(\tilde\mu^o,\tilde\mu,\tilde\eta)$ in $C^{1,\alpha}(\partial\Omega^o)\times C^{0,\alpha}(\partial\Omega^i)\times C^{0,\alpha}(\partial\Omega^i)$.  
\end{prop}
\proof 
Let $T$ be as in \eqref{T}.  By Proposition \ref{C} there exists $(\tilde\mu^o,\tilde\mu,\tilde\eta)$ in $C^{1,\alpha}(\partial\Omega^o)\times C^{0}(\partial\Omega^i)\times C^{0}(\partial\Omega^i)$ such that $(\tilde\mu^o,\tilde\mu,\tilde\eta)=T(\tilde\mu^o,\tilde\mu,\tilde\eta)$. Then, by the mapping properties of integral operators with real analytic kernels and no singularities we have that $w^-_{\Omega^i}[\tilde\mu]_{|\partial\Omega^o}$ and $v_{\Omega^i}[\tilde\eta]_{|\partial\Omega^o}$ belong to $C^{1,\alpha}(\partial\Omega^o)$, that  $w^+_{\Omega^o}[\tilde\mu^o]_{|\partial\Omega^i}$ belongs to $C^{1,\alpha}(\partial\Omega^i)$, and that $\nu_{\Omega^i}\cdot\nabla w^+_{\Omega^o}[\psi]_{|\partial\Omega^i}$ belongs to  $C^{0,\alpha}(\partial\Omega^i)$ for all $\psi\in C^0(\partial\Omega^o)$. By a classical result in potential theory (cf., {\it e.g.}, Miranda \cite[Chap.~II, \S14, III]{Mi70}) we have that $v_{\Omega^i}[\tilde\eta]_{|\partial\Omega^i}\in C^{0,\alpha}(\partial\Omega^i)$ and by Lemma \ref{schauder} we have that $W_{\Omega^i}[\mu]\in C^{0,\alpha}(\partial\Omega^i)$. Then, by the invertibility of $\frac{1}{2}I_{\Omega^o}+W_{\Omega^o}$ in $C^{1,\alpha}(\partial\Omega^o)$ and of $\frac{1}{2}I_{\Omega^i}+W_{\Omega^i}$ in $C^{0,\alpha}(\partial\Omega^i)$ (cf.~Lemma \ref{schauder} and \ref{ker}),  by the invertibility of $J$ in $C^{0,\alpha}(\partial\Omega^i)$ (cf.~Lemma \ref{J}), and by assumption \eqref{assumption2} it follows that $T(\tilde\mu^o,\tilde\mu,\tilde\eta)\in C^{1,\alpha}(\partial\Omega^o)\times C^{0,\alpha}(\partial\Omega^i)\times C^{0,\alpha}(\partial\Omega^i)$. Thus $(\tilde\mu^o,\tilde\mu,\tilde\eta)\in C^{1,\alpha}(\partial\Omega^o)\times C^{0,\alpha}(\partial\Omega^i)\times C^{0,\alpha}(\partial\Omega^i)$ and our proof is complete.\qed

\vspace{\baselineskip}

In the following Theorem \ref{existence} we show that under conditions \eqref{assumption0} and \eqref{assumption1} there exists a pair of functions $(\tilde u^o,\tilde u^i)\in C^0(\mathrm{cl}\Omega^o\setminus\Omega^i)\times C^0(\mathrm{cl}\Omega^i)$ which satisfy the first four conditions of problem \eqref{nlprob}  in the classical sense and which satisfies the fifth condition of \eqref{nlprob} in a certain weak sense which we now specify.  To do so, we  define  the distribution $[\nu_{\Omega^i}\cdot\nabla\tilde w^o- \nu_{\Omega^i}\cdot\nabla\tilde w^i]_w$ for all pair of functions $(\tilde w^o, \tilde w^i)\in C^0_{\mathrm{harm}}(\mathrm{cl}\Omega^o\setminus\Omega^i)\times C^0_{\mathrm{harm}}(\mathrm{cl}\Omega^i)$ (see also definition \eqref{Charm}).

\begin{defn}\label{weak}
Let $(\tilde w^o, \tilde w^i)$ be a pair of functions of $C^0_{\mathrm{harm}}(\mathrm{cl}\Omega^o\setminus\Omega^i)\times C^0_{\mathrm{harm}}(\mathrm{cl}\Omega^i)$. Then  $[\nu_{\Omega^i}\cdot\nabla\tilde w^o- \nu_{\Omega^i}\cdot\nabla\tilde w^i]_w$ denotes the distribution on $\Omega^o$ defined by
\[
\begin{split}
&\langle [\nu_{\Omega^i}\cdot\nabla\tilde w^o- \nu_{\Omega^i}\cdot\nabla\tilde w^i]_w, \phi \rangle\\
&\quad\equiv\int_{\partial\Omega^i}(\tilde w^o_{|\partial\Omega^i}-\tilde w^i_{|\partial\Omega^i})(\nu_{\Omega^i}\cdot\nabla \phi_{|\partial\Omega^i})\, d\sigma+\int_{\Omega^o\setminus\Omega^i}\tilde w^o\ \Delta \phi\, dx+\int_{\Omega^i}\tilde w^i\ \Delta \phi\, dx
\end{split}
\] 
for all test functions $\phi\in C^\infty_c(\Omega^o)$.
\end{defn}

One immediately verifies that the map which takes $(\tilde w^o, \tilde w^i)$ to $[\nu_{\Omega^i}\cdot\nabla\tilde w^o- \nu_{\Omega^i}\cdot\nabla\tilde w^i]_w$ is continuous. Namely we have the following.
\begin{lem}\label{continuity} Let $(\tilde w^o, \tilde w^i)$ be a pair of functions of $C^0_{\mathrm{harm}}(\mathrm{cl}\Omega^o\setminus\Omega^i)\times C^0_{\mathrm{harm}}(\mathrm{cl}\Omega^i)$ and let $\{(\tilde w^o_j, \tilde w^i_j)\}_{j\in\mathbb{N}}$ be a sequence in $C^0_{\mathrm{harm}}(\mathrm{cl}\Omega^o\setminus\Omega^i)\times C^0_{\mathrm{harm}}(\mathrm{cl}\Omega^i)$ such that $\lim_{j\to+\infty}\tilde w^o_j=\tilde w^o$ in $C^0(\mathrm{cl}\Omega^o\setminus\Omega^i)$ and  $\lim_{j\to+\infty}\tilde w^i_j=\tilde w^i$ in $C^0(\mathrm{cl}\Omega^i)$. Then  
\[
\lim_{j\to+\infty}\langle [\nu_{\Omega^i}\cdot\nabla\tilde w^o_j- \nu_{\Omega^i}\cdot\nabla\tilde w^i_j]_w,\phi\rangle=\langle [\nu_{\Omega^i}\cdot\nabla\tilde w^o- \nu_{\Omega^i}\cdot\nabla\tilde w^i]_w,\phi\rangle\qquad\forall \phi\in C^\infty_c(\Omega^o)\,.
\]
\end{lem}
Moreover, if $(w^o,  w^i)$ belongs to $C^1_{\mathrm{harm}}(\mathrm{cl}\Omega^o\setminus\Omega^i)\times C^1_{\mathrm{harm}}(\mathrm{cl}\Omega^i)$, then  $[\nu_{\Omega^i}\cdot\nabla w^o- \nu_{\Omega^i}\cdot\nabla w^i]_w$ coincides with $(\nu_{\Omega^i}\cdot\nabla  w^o-\nu_{\Omega^i}\cdot\nabla   w^i)_{|\partial\Omega^i}$. Namely we have
\[
\langle [\nu_{\Omega^i}\cdot\nabla  w^o-  \nu_{\Omega^i}\cdot\nabla w^i]_w, \phi \rangle=\int_{\partial\Omega^i}\left(\nu_{\Omega^i}\cdot\nabla w^o(x)- \nu_{\Omega^i}\cdot\nabla w^i(x)\right)\phi(x)\,d\sigma_x\
\]
for all $\phi\in C^\infty_c(\Omega^o)$ and for all pair of functions $( w^o,  w^i)\in C^1_{\mathrm{harm}}(\mathrm{cl}\Omega^o\setminus\Omega^i)\times C^1_{\mathrm{harm}}(\mathrm{cl}\Omega^i)$. Then we can prove that $[\nu_{\Omega^i}\cdot\nabla\tilde w^o- \nu_{\Omega^i}\cdot\nabla\tilde w^i]_w$ is supported on $\partial\Omega^i$.

\begin{lem}\label{bounadary}
For all $(\tilde w^o, \tilde w^i)\in C^0_{\mathrm{harm}}(\mathrm{cl}\Omega^o\setminus\Omega^i)\times C^0_{\mathrm{harm}}(\mathrm{cl}\Omega^i)$  the support of $[\nu_{\Omega^i}\cdot\nabla\tilde w^o- \nu_{\Omega^i}\cdot\nabla\tilde w^i]_w$ is contained in $\partial\Omega^i$. 
\end{lem}
\proof
By a classical argument  one can prove that there exists a sequence $\{(w_j^o,w_j^i)\}_{j\in\mathbb{N}}$ in $C^{1,\alpha}_{\mathrm{harm}}(\mathrm{cl}\Omega^o\setminus\Omega^i)\times C^{1,\alpha}_{\mathrm{harm}}(\mathrm{cl}\Omega^i)$ such that  $\lim_{j\to+\infty}w_j^o=\tilde w^o$ in $C^{0}(\mathrm{cl}\Omega^o\setminus\Omega^i)$ and $\lim_{j\to+\infty}w_j^i=\tilde w^i$ in $C^{0}(\mathrm{cl}\Omega^i)$. Let $\phi_0\in C^\infty_c(\Omega^o)$ be such that $\phi_{0|\partial\Omega^i}=0$.  Then we have
\[
\langle [\nu_{\Omega^i}\cdot\nabla w^o_j-  \nu_{\Omega^i}\cdot\nabla  w^i_j]_w, \phi_0 \rangle=\int_{\partial\Omega^i}\left(\nu_{\Omega^i}\cdot\nabla w^o_j(x)- \nu_{\Omega^i}\cdot\nabla w^i_j(x)\right)\phi_0(x)\,d\sigma_x=0
\] 
for all $j\in\mathbb{N}$. Moreover $\lim_{j\to\infty}\langle [\nu_{\Omega^i}\cdot\nabla w^o_j-  \nu_{\Omega^i}\cdot\nabla  w^i_j]_w, \phi_0 \rangle=\langle [\nu_{\Omega^i}\cdot\nabla \tilde w^o-  \nu_{\Omega^i}\cdot\nabla \tilde w^i]_w, \phi_0 \rangle$ by Lemma \ref{continuity}, and thus  $\langle[\nu_{\Omega^i}\cdot\nabla\tilde w^o- \nu_{\Omega^i}\cdot\nabla\tilde w^i]_w, \phi_0 \rangle=0$.
\qed

\vspace{\baselineskip}

We are now ready to prove the main result of this section.

\begin{thm}\label{existence}
Assume that $F$ and $G$ satisfy \eqref{assumption0} and  \eqref{assumption1}. Then there exists $(\tilde{u}^o,\tilde{u}^i)\in C^0(\mathrm{cl}\Omega^o\setminus\Omega^i)\times C^0(\mathrm{cl}\Omega^i)$ such that 
\begin{equation}\label{existence.eq1}
\left\{
\begin{array}{ll}
\Delta \tilde u^o=0&\text{in }\Omega^o\setminus\mathrm{cl}\Omega^i\,,\\
\Delta \tilde u^i=0&\text{in }\Omega^o\,,\\
\tilde u^o(x)=f^o(x)&\text{for all }x\in\partial\Omega^o\,,\\
\tilde u^o(x)=F(x,\tilde u^i(x))&\text{for all }x\in\partial\Omega^i\,,\\
\langle [\nu_{\Omega^i}\cdot\nabla\tilde u^o-\nu_{\Omega^i}\cdot\nabla\tilde u^i)]_w,\phi\rangle=\int_{\partial\Omega^i} G (x,\tilde u^i(x))\phi(x)\,d\sigma_x&\text{for all }\phi\in C^{\infty}_c(\Omega^o)\,.
\end{array}
\right.
\end{equation}  
\end{thm}

\proof Let $(\tilde\mu^o,\tilde\mu,\tilde\eta)\in C^{1,\alpha}(\partial\Omega^o)\times C^{0}(\partial\Omega^i)\times C^{0}(\partial\Omega^i)$ be as in Proposition \ref{existenceC} and define 
\[
\tilde{u}^o\equiv(w^+_{\Omega^o}[\tilde\mu^o]+w^-_{\Omega^i}[\tilde\mu]+v^-_{\Omega^i}[\tilde\eta])_{|\mathrm{cl}\Omega^o\setminus\Omega^i}\,,\quad \tilde{u}^i\equiv w^+_{\Omega^i}[\tilde\mu]\,.
\]
Then the pair $(\tilde{u}^o,\tilde{u}^i)$ belongs to $C^0(\mathrm{cl}\Omega^o\setminus\Omega^i)\times C^0(\mathrm{cl}\Omega^i)$ (cf. Folland \cite[Chap.~3]{Fo95}) and satisfies the first four conditions of \eqref{existence.eq1} (see also \eqref{jump}). We now prove that $(\tilde{u}^o,\tilde{u}^i)$ satisfies also the fifth condition of \eqref{existence.eq1}.

 By a standard argument  one proves that there exists a sequence $\{v_j^i\}_{j\in\mathbb{N}}$    in $C^{1,\alpha}_{\mathrm{harm}}(\mathrm{cl}\Omega^i)$ such that  
 \begin{equation}\label{existence.eq5}
 \lim_{j\to+\infty}v^i_j=\tilde{u}^i\quad\text{in }C^0(\mathrm{cl}\Omega^i)\,.
 \end{equation}
 By Lemmas \ref{fredC} and \ref{ker}  we have that $\frac{1}{2}I_{\Omega^i}+W_{\Omega^i}$  is an isomorphism from  $C^0(\partial\Omega^i)$ to itself and from $C^{1,\alpha}(\partial\Omega^i)$ to itself. Then, by \eqref{jump}  one verifies that there exists  $\mu_j\in C^{1,\alpha}(\partial\Omega^i)$ such that $v^i_j= w^+_{\Omega^i}[\mu_j]$ for all $j\in\mathbb{N}$. Moreover, by the continuity of $(\frac{1}{2}I_{\Omega^i}+W_{\Omega^i})^{(-1)}$ from $C^0(\partial\Omega^i)$ to itself, we have
\begin{equation}\label{existence.eq2}
\lim_{j\to +\infty}\mu_j=\lim_{j\to +\infty}(\frac{1}{2}I_{\Omega^i}+W_{\Omega^i})^{(-1)}v^i_{j|\partial\Omega^i}=(\frac{1}{2}I_{\Omega^i}+W_{\Omega^i})^{(-1)}\tilde{u}^i_{|\partial\Omega^i}=\tilde\mu\quad \text{in }C^0(\partial\Omega^i)\,.
\end{equation}
Then we set 
\[
\mu^o_j\equiv(\frac{1}{2}I_{\Omega^o}+W_{\Omega^o})^{(-1)}(f^o-w^-_{\Omega^i}[\mu_j]_{|\partial\Omega^o}-v_{\Omega^i}[\tilde\eta]_{|\partial\Omega^o})\quad\forall j\in\mathbb{N}\,.
\]  By Lemmas \ref{fredC} and \ref{ker},  we have that $\frac{1}{2}I_{\Omega^o}+W_{\Omega^o}$  is an isomorphism from $C^{1,\alpha}(\partial\Omega^o)$ to itself and from $C^{0}(\partial\Omega^o)$ to itself. In particular, $(\frac{1}{2}I_{\Omega^o}+W_{\Omega^o})^{(-1)}$ is continuous from  $C^{0}(\partial\Omega^o)$ to itself and maps $C^{1,\alpha}(\partial\Omega^o)$ to itself. Moreover, by the properties of integral operators with real analytic kernel and no singularities $w^-_{\Omega^i}[\cdot]_{|\partial\Omega^o}$ is continuous from $C^{0}(\partial\Omega^i)$ to $C^{1,\alpha}(\partial\Omega^o)$ and $f^o-v_{\Omega^i}[\tilde\eta]_{|\partial\Omega^o}$ belongs to $C^{1,\alpha}(\partial\Omega^o)$. It follows that $\mu^o_j\in C^{1,\alpha}(\partial\Omega^o)$ for all $j\in\mathbb{N}$ and that
\begin{equation}\label{existence.eq3}
\begin{split}
\lim_{j\to +\infty}\mu^o_j&=\lim_{j\to +\infty}(\frac{1}{2}I_{\Omega^o}+W_{\Omega^o})^{(-1)}(f^o-w^-_{\Omega^i}[\mu_j]_{|\partial\Omega^o}-v_{\Omega^i}[\tilde\eta]_{|\partial\Omega^o})\\
&=(\frac{1}{2}I_{\Omega^o}+W_{\Omega^o})^{(-1)}(f^o-w^-_{\Omega^i}[\tilde\mu]_{|\partial\Omega^o}-v_{\Omega^i}[\tilde\eta]_{|\partial\Omega^o})=\tilde\mu^o\qquad \text{in }C^{0}(\partial\Omega^o)\,.
\end{split}
\end{equation}
Now let
\[
v^o_j\equiv (w^+_{\Omega^o}[\mu^o_j]+w^-_{\Omega^i}[\mu_j]+v^-_{\Omega^i}[\tilde \eta])_{|\mathrm{cl}\Omega^o\setminus\Omega^i}\quad\forall j\in\mathbb{N}\,.
\]
By classical potential theory  $v^o_j\in C^{1,\alpha}(\mathrm{cl}\Omega^o\setminus\Omega^i)$  (cf., {\it e.g.}, Miranda \cite{Mi65}). Moreover, by \eqref{jump} we have
\[
v^o_{j|\partial\Omega^o}=(\frac{1}{2}I_{\Omega^o}+W_{\Omega^o})\mu^o_j+w^-_{\Omega^i}[\mu_j]_{|\partial\Omega^o}+v^-_{\Omega^i}[\tilde \eta]_{|\partial\Omega^o}
\]
and
\[
v^o_{j|\partial\Omega^i}=w^+_{\Omega^o}[\mu^o_j]_{|\partial\Omega^i}+(-\frac{1}{2}I_{\Omega^i}+W_{\Omega^i})\mu_j+v^-_{\Omega^i}[\tilde \eta]_{|\partial\Omega^i}\,.
\]
Then, by \eqref{existence.eq2} and \eqref{existence.eq3}, by the continuity of  $\frac{1}{2}I_{\Omega^o}+W_{\Omega^o}$ from $C^0(\partial\Omega^o)$ to itself and of $-\frac{1}{2}I_{\Omega^o}+W_{\Omega^o}$ from $C^0(\partial\Omega^i)$ to itself (cf.~Lemma \ref{schauder}), and by the properties of integral operators with real analytic kernels and no singularity, we deduce that 
\[
\lim_{j\to +\infty}v^o_{j|\partial\Omega^o}=\tilde u^o_{|\partial\Omega^o}\quad\text{in }C^{0}(\partial\Omega^o)
\] 
and
\[
\lim_{j\to +\infty}v^o_{j|\partial\Omega^i}=\tilde u^o_{|\partial\Omega^i}\quad\text{in }C^{0}(\partial\Omega^i)\,.
\] 
It follows that 
\begin{equation}\label{existence.eq4}
\lim_{j\to +\infty}v^o_j=\tilde u^o\quad\text{in }C^{0}(\mathrm{cl}\Omega^o\setminus\Omega^i)\,.
\end{equation}
In addition, by the jump formulas  \eqref{jump} and \eqref{nojump} and by the validity of equality \eqref{integral} for $(\mu^o,\mu,\eta)=(\tilde\mu^o,\tilde\mu,\tilde\eta)$ (cf.~Proposition \ref{existenceC}) one verifies that the pair $(v_j^o,v_j^i)$ satisfies the equality 
\[
\nu_{\Omega^i}\cdot\nabla v^o_{j}(x)-\nu_{\Omega^i}\cdot\nabla v^i_j(x)=G(x,\tilde{u}^i(x))+\nu_{\Omega^i}(x)\cdot\nabla w^+_{\Omega^o}[\mu^o_j-\tilde\mu^o](x)\quad\forall x\in\partial\Omega^i 
\] for all $j\in\mathbb{N}$. Hence, by the continuity of the map  from $C^0(\partial\Omega^o)$ to $C^0(\partial\Omega^i)$ which takes $\phi$ to $\nu_{\Omega^i}\cdot\nabla w^+_{\Omega^o}[\phi]_{|\partial\Omega^i}$ and by the limit relation in \eqref{existence.eq3} we have that 
\[
\lim_{j\to +\infty}(\nu_{\Omega^i}\cdot\nabla v^o_{j}-\nu_{\Omega^i}\cdot\nabla v^i_j)_{|\partial\Omega^i}=\mathcal{F}_G\tilde{u}^i\qquad\text{in }C^{0}(\partial\Omega^i)\,.
\] 
Thus,  by Lemma \ref{continuity}, by the limit relations in \eqref{existence.eq5} and \eqref{existence.eq4}, and by the membership of $(v^o_j,v^i_j)$ in $C^{1,\alpha}_{\mathrm{harm}}(\mathrm{cl}\Omega^o\setminus\Omega^i)\times C^{1,\alpha}_{\mathrm{harm}}(\mathrm{cl}\Omega^i)$ for all $j\in\mathbb{N}$, it follows that $(\tilde u^o,\tilde u^i)$ satisfies the fifth condition in problem \eqref{existence.eq1}. The theorem is now proved. \qed

\vspace{\baselineskip}

If in addition $F$ and $G$ satisfy assumption \eqref{assumption2}, then the pair $(\tilde{u}^o,\tilde{u}^i)$ belongs to $C^{0,\alpha}(\mathrm{cl}\Omega^o\setminus\Omega^i)\times C^{0,\alpha}(\mathrm{cl}\Omega^i)$.

\begin{thm}\label{existenceC0a}
Assume that $F$ and $G$ satisfy \eqref{assumption0}, \eqref{assumption1}, and \eqref{assumption2}. Then there exists $(\tilde{u}^o,\tilde{u}^i)\in C^{0,\alpha}(\mathrm{cl}\Omega^o\setminus\Omega^i)\times C^{0,\alpha}(\mathrm{cl}\Omega^i)$ which satisfy the conditions in \eqref{existence.eq1}. \end{thm}

\proof If $(\tilde\mu^o,\tilde\mu,\tilde\eta)\in C^{1,\alpha}(\partial\Omega^o)\times C^{0,\alpha}(\partial\Omega^i)\times C^{0,\alpha}(\partial\Omega^i)$ is as in Proposition \ref{existenceC0a} and 
\[
\tilde{u}^o\equiv(w^+_{\Omega^o}[\tilde\mu^o]+w^-_{\Omega^i}[\tilde\mu]+v^-_{\Omega^i}[\tilde\eta])_{|\mathrm{cl}\Omega^o\setminus\Omega^i}\,,\quad \tilde{u}^i\equiv w^+_{\Omega^i}[\tilde\mu]\,,
\]
then the pair $(\tilde{u}^o,\tilde{u}^i)$ belongs to $C^{0,\alpha}(\mathrm{cl}\Omega^o\setminus\Omega^i)\times C^{0,\alpha}(\mathrm{cl}\Omega^i)$ (cf.~Miranda \cite{Mi65}) and we can prove that it satisfies the  conditions of \eqref{existence.eq1} by arguing as in the proof of Theorem \ref{existence}.
\qed

\section{Existence result for problem \eqref{nlprob.small}}\label{sec.small}

We now  fix a real number $\lambda>0$ and a continuous function $\Phi$  from $\partial\Omega^i\times\mathbb{R}$ to $\mathbb{R}$. Then  we assume that $F=\lambda\mathrm{id}_\mathbb{R}+\epsilon\Phi$, where $\epsilon$ is a multiplicative real parameter. Our aim is to study the nonlinear transmission problem \eqref{nlprob.small} for $\epsilon$ small. To do so, we find convenient to introduce the following technical assumption:
\begin{equation}\label{ass.omega}
\begin{split}
&\text{the map from $C^{0,\alpha}(\partial\Omega^i)$ to $C^{1,\alpha}(\partial\Omega^i)$ which}\\
&\text{takes $\eta$ to $v_{\Omega^i}[\eta]_{|\partial\Omega^i}$ is an isomorphism.}
\end{split}
\end{equation}
We observe that assumption \eqref{ass.omega} holds for all domains $\Omega^i$ in $\mathbb{R}^n$ if $n\ge 3$, and does not old in $\mathbb{R}^2$ only in exceptional cases. Indeed we have the following classical result.

\begin{prop}
Let $\Omega$ be an open bounded connected subset of $\mathbb{R}^2$ of class $C^{1,\alpha}$. Then there exists  a unique function $\psi_\Omega$ in $C^{0,\alpha}(\partial\Omega)$ such that 
\[
\text{$v_{\Omega}[\psi_\Omega]_{|\partial\Omega}$ is constant and $\int_{\partial\Omega}\psi_\Omega\,d\sigma=1$.}
\]  
Moreover, the following  statements hold.
\begin{enumerate}
\item[(i)] If $v_{\Omega}[\psi_\Omega]_{|\partial\Omega}\neq 0$, then $v_{\Omega}[\cdot]_{|\partial\Omega}$ is an isomorphism from $C^{0,\alpha}(\partial\Omega)$ to $C^{1,\alpha}(\partial\Omega)$.
\item[(ii)] If  $v_{\Omega}[\psi_\Omega]_{|\partial\Omega}=0$ and $r\in]0,+\infty[\setminus\{1\}$, then  $v_{r\Omega}[\cdot]_{|r\partial\Omega}$ is an isomorphism from $C^{0,\alpha}(r\partial\Omega)$ to $C^{1,\alpha}(r\partial\Omega)$.
\end{enumerate}
 \end{prop}

In Lemma \ref{representation.small} here below we introduce an isomorphism between $C^{1,\alpha}(\partial\Omega^o)\times C^{0,\alpha}(\partial\Omega^i)\times C^{0,\alpha}(\partial\Omega^i)$ and $C^{1,\alpha}_\mathrm{harm}(\mathrm{cl}\Omega^o\setminus\Omega^i)\times C^{1,\alpha}_\mathrm{harm}(\mathrm{cl}\Omega^i)$ (cf.~definition \eqref{Charm}).

\begin{lem}\label{representation.small}
Let $U\equiv (U^o,U^i)$ denote  the operator from $C^{1,\alpha}(\partial\Omega^o)\times C^{0,\alpha}(\partial\Omega^i)\times C^{0,\alpha}(\partial\Omega^i)$ to $C^{1,\alpha}_\mathrm{harm}(\mathrm{cl}\Omega^o\setminus\Omega^i)\times C^{1,\alpha}_\mathrm{harm}(\mathrm{cl}\Omega^i)$ which takes   $(\mu^o,\eta^o,\eta^i)$ to the pair $(U^o[\mu^o,\eta^o,\eta^i], U^i[\mu^o,\eta^o,\eta^i])$ given by 
\[
 U^o[\mu^o,\eta^o,\eta^i]\equiv(w^+_{\Omega^o}[\mu^o]+v^-_{\Omega^i}[\eta^o])_{|\mathrm{cl}\Omega^o\setminus\Omega^i}\,,\quad
U^i[\mu^o,\eta^o,\eta^i]\equiv \lambda^{-1}w^+_{\Omega^o}[\mu^o]_{|\mathrm{cl}\Omega^i}+\lambda^{-1}v^+_{\Omega^i}[\eta^i]\,.
\] 
Then $U$ is a linear isomorphism.
\end{lem}
\proof By the mapping properties of the single and double layer potentials one verifies that the operator $U$  is  continuous  from $C^{1,\alpha}(\partial\Omega^o)\times C^{0,\alpha}(\partial\Omega^i)\times C^{0,\alpha}(\partial\Omega^i)$ to $C^{1,\alpha}_\mathrm{harm}(\mathrm{cl}\Omega^o\setminus\Omega^i)\times C^{1,\alpha}_\mathrm{harm}(\mathrm{cl}\Omega^i)$ (cf.~Section \ref{sec.preliminaries}, see also Miranda \cite{Mi65}).

Therefore, if we prove that $U$ is  one-to-one and onto, we can deduce by the open mapping theorem  that $U$ is an isomorphism from $C^{1,\alpha}(\partial\Omega^o)\times C^{0,\alpha}(\partial\Omega^i)\times C^{0,\alpha}(\partial\Omega^i)$ to $C^{1,\alpha}_\mathrm{harm}(\mathrm{cl}\Omega^o\setminus\Omega^i)\times C^{1,\alpha}_\mathrm{harm}(\mathrm{cl}\Omega^i)$. So let $(\phi^o,\phi^i)\in C^{1,\alpha}_\mathrm{harm}(\mathrm{cl}\Omega^o\setminus\Omega^i)\times C^{1,\alpha}_\mathrm{harm}(\mathrm{cl}\Omega^i)$. We show that there exists unique  triple $(\mu^o,\eta^o,\eta^i)\in C^{1,\alpha}(\partial\Omega^o)\times C^{0,\alpha}(\partial\Omega^i)\times C^{0,\alpha}(\partial\Omega^i)$ such that $(U^o[\mu^o,\eta^o,\eta^i], U^i[\mu^o,\eta^o,\eta^i])=(\phi^o,\phi^i)$.   We first  consider $U^o[\mu^o,\eta^o,\eta^i]=\phi^o$ and we verify that there exists unique $(\mu^o,\eta^o)$ such that 
\begin{equation}\label{representation.small.eq1}
(w^+_{\Omega^o}[\mu^o]+v^-_{\Omega^i}[\eta^o])_{|\mathrm{cl}\Omega^o\setminus\Omega^i}=\phi^o\,.
\end{equation} 
By the uniqueness of the classical solution of the Neumann-Dirichlet mixed problem and by the jump properties of the single and double layer potentials (cf.~equality \eqref{jump}), equation \eqref{representation.small.eq1} is equivalent to 
\[
\begin{split}
&(\frac{1}{2}I_{\Omega^o}+W_{\Omega^o})\mu^o+v^-_{\Omega^i}[\eta^o]_{|\partial\Omega^o}=\phi^o_{|\partial\Omega^o}\,,\\
&(\frac{1}{2}I_{\Omega^i}+W^*_{\Omega^i})\eta^o+\nu_{\Omega^i}\cdot\nabla w^+_{\Omega^o}[\mu^o]_{|\partial\Omega^i}=\nu_{\Omega^i}\cdot\nabla \phi^o_{|\partial\Omega^i}.
\end{split}
\]
 By classical potential theory, the operator $(\frac{1}{2}I_{\Omega^o}+W_{\Omega^o})$ is an isomorphism from $C^{1,\alpha}(\partial\Omega^o)$ to itself and  the operator $(\frac{1}{2}I_{\Omega^i}+W^*_{\Omega^i})$ is an isomorphism from $C^{0,\alpha}(\partial\Omega^i)$ to itself (cf.~Lemmas \ref{fredC} and \ref{ker}). Then, by the compactness properties of the integral operators with real analytic kernel and with no singularities and by standard properties of Fredholm operators, we deduce that the map which takes $(\mu,\eta)$ to $((\frac{1}{2}I_{\Omega^o}+W_{\Omega^o})\mu+v^-_{\Omega^i}[\eta]_{|\partial\Omega^o}, (\frac{1}{2}I_{\Omega^i}+W^*_{\Omega^i})\eta+\nu_{\Omega^i}\cdot\nabla w^+_{\Omega^o}[\mu^o]_{|\partial\Omega^i})$ is a Fredholm operator of index $0$ from $C^{1,\alpha}(\partial\Omega^o)\times C^{0,\alpha}(\partial\Omega^i)$ to itself. Thus,  to prove the existence and uniqueness of $(\mu^o,\eta^o)\in C^{1,\alpha}(\partial\Omega^o)\times  C^{0,\alpha}(\partial\Omega^i)$ which satisfies \eqref{representation.small.eq1} it suffices to show that $(\mu,\eta)=(0,0)$ when 
 \begin{equation}\label{representation.small.eq2}
\Bigl((\frac{1}{2}I_{\Omega^o}+W_{\Omega^o})\mu+v^-_{\Omega^i}[\eta]_{|\partial\Omega^o}\,,\; (\frac{1}{2}I_{\Omega^i}+W^*_{\Omega^i})\eta+\nu_{\Omega^i}\cdot\nabla w^+_{\Omega^o}[\mu]_{|\partial\Omega^i}\Bigr)=(0,0)\,.
 \end{equation}
 If \eqref{representation.small.eq2} holds,  then   $(w^+_{\Omega^o}[\mu]+v^-_{\Omega^i}[\eta])_{|\mathrm{cl}\Omega^o\setminus\Omega^i}=0$ by the uniqueness of the classical solution of the Neumann-Dirichlet mixed problem and by \eqref{jump}. Hence $w^+_{\Omega^o}[\mu]+v^+ _{\Omega^i}[\eta]=0$ in $\mathrm{cl}\Omega^i$ by the uniqueness of the classical solution of the  Dirichlet  problem in $\Omega^i$ and by the continuity of $(w^+_{\Omega^o}[\mu]+v_{\Omega^i}[\eta])_{|\mathrm{cl}\Omega^o}$ (cf.~Section \ref{sec.preliminaries}). Then, by the jump properties of the single layer potential (cf.~equality \eqref{jump})  we have that 
 \[
 \begin{split}
 \eta&=\nu_{\Omega^i}\cdot\nabla v^-_{\Omega^i}[\eta]_{\partial\Omega^i}-\nu_{\Omega^i}\cdot\nabla v^+_{\Omega^i}[\eta]_{|\partial\Omega^i}\\
 &=\nu_{\Omega^i}\cdot\nabla(w^+_{\Omega^o}[\mu]+v^-_{\Omega^i}[\eta])_{|\partial\Omega^i}-\nu_{\Omega^i}\cdot\nabla(w^+_{\Omega^o}[\mu]+v^+_{\Omega^i}[\eta])_{|\partial\Omega^i}=0
 \end{split}
 \]
 and  thus $\mu=0$ by the first equality in \eqref{representation.small.eq2} and by Lemma \ref{ker} (iii). Now, to complete the proof we observe that $U^i[\mu^o,\eta^o,\eta^i]=\phi^i$ is equivalent to $v^+_{\Omega^i}[\eta^i]=\lambda\phi^i-w^+_{\Omega^o}[\mu^o]_{|\mathrm{cl}\Omega^i}$ and the existence and uniqueness of $\eta^i$ is guaranteed by the assumption in \eqref{ass.omega}.
\qed

\vspace{\baselineskip}

In the following Lemma \ref{Jl} we introduce an auxiliary operator which we denote by $J_\lambda$.

\begin{lem}\label{Jl}
Let 
\[
J_\lambda[\eta]\equiv\left(\frac{1}{2}I_{\Omega^i}+\frac{\lambda-1}{\lambda+1}W^*_{\Omega^i}\right)\eta-\frac{\lambda-1}{\lambda+1}\nu_{\Omega^i}\cdot\nabla w^+_{\Omega^o}\biggl[(\frac{1}{2}I_{\Omega^o}+W_{\Omega^o})^{(-1)}v_{\Omega^i}[\eta]_{|\partial\Omega^o}\biggr]_{|\partial\Omega^i}
\] for all $\eta\in L^2(\partial\Omega^i)$. Then the map which takes $\eta$ to $J_\lambda[\eta]$ is an isomorphism from $L^2(\partial\Omega^i)$ to itself, from $C^0(\partial\Omega^i)$ to itself, and from $C^{0,\alpha}(\partial\Omega^i)$ to itself.
\end{lem}
\proof By the properties of integral operators with real analytic kernels and no singularity, by the invertibility of $\frac{1}{2}I_{\Omega^o}+W_{\Omega^o}$ in $C^{1,\alpha}(\partial\Omega^o)$ (cf.~Lemma \ref{fredC} and Lemma \ref{ker}), and by the continuity of the map $w^+_{\Omega^o}[\cdot]$ from $C^{1,\alpha}(\partial\Omega^o)$ to $C^{1,\alpha}(\mathrm{cl}\Omega^o)$ (cf., {\it e.g.}, Miranda \cite{Mi65}), one deduces that the operator which takes $\eta$ to  
\begin{equation}\label{Jl.eq1}
\nu_{\Omega^i}\cdot\nabla w^+_{\Omega^o}\biggl[(\frac{1}{2}I_{\Omega^o}+W_{\Omega^o})^{(-1)}v_{\Omega^i}[\eta]_{|\partial\Omega^o}\biggr]_{|\partial\Omega^i}
\end{equation}
is continuous from $L^2(\partial\Omega^i)$ to $C^{0,\alpha}(\partial\Omega^i)$. Then, by the invertibility of $\frac{1}{2}I_{\Omega^i}+\frac{\lambda-1}{\lambda+1}W^*_{\Omega^i}$ in $L^2(\partial\Omega^i)$ (cf.~Lemma \ref{ilw}) it follows that $J_\lambda$ is a Fredholm operator of index $0$ from $L^2(\partial\Omega^i)$ to itself. Thus, to show that $J_\lambda$ is invertible from $L^2(\partial\Omega^i)$ to itself it suffices to prove that $J_\lambda[\eta]=0$ implies $\eta=0$. Now, if $\eta\in L^2(\partial\Omega^i)$ and $J_\lambda[\eta]=0$, then $(\frac{1}{2}I_{\Omega^i}+\frac{\lambda-1}{\lambda+1}W^*_{\Omega^i})\eta\in C^{0,\alpha}(\partial\Omega^i)$ by the membership of \eqref{Jl.eq1} in $C^{0,\alpha}(\partial\Omega^i)$, and thus $\eta\in C^{0,\alpha}(\partial\Omega^i)$ by Lemma \ref{regularity}.  Then, by taking $\mu^o\equiv -(\frac{1}{2}I_{\Omega^o}+W_{\Omega^o})^{(-1)}v_{\Omega^i}[\eta]_{|\partial\Omega^o}$ and by a straightforward calculation based on \eqref{jump} one verifies that 
$U^o[\mu^o,\eta,\eta]_{|\partial\Omega^o}=0$, $U^o[\mu^o,\eta,\eta]_{|\partial\Omega^i}=\lambda U^i[\mu^o,\eta,\eta]_{|\partial\Omega^i}$,
and $\nu_{\Omega^i}\cdot\nabla U^o[\mu^o,\eta,\eta]_{|\partial\Omega^i}= \nu_{\Omega^i}\cdot\nabla U^i[\mu^o,\eta,\eta]_{|\partial\Omega^i}$
(where  $U^o[\mu^o,\eta,\eta]$ and $U^i[\mu^o,\eta,\eta]$ are defined as in Lemma \ref{representation.small}). Then, by the uniqueness of the solution of the linear perfect contact problem we have $U^o[\mu^o,\eta,\eta]=0$ and $U^i[\mu^o,\eta,\eta]=0$. Accordingly, Lemma \ref{representation.small} implies that $\eta=0$. 

To prove that $J_\lambda$ is invertible from $C^0(\partial\Omega^i)$ to itself, we first observe that $J_\lambda$ is continuous from $C^0(\partial\Omega^i)$ to itself (see Lemma \ref{ilw}). Moreover, if $\eta\in L^2(\partial\Omega^i)$ and $J_\lambda[\eta]\in C^0(\partial\Omega^i)$ then $(\frac{1}{2}I_{\Omega^i}+\frac{\lambda-1}{\lambda+1}W^*_{\Omega^i})\eta\in C^{0}(\partial\Omega^i)$ by the membership of \eqref{Jl.eq1} in $C^{0,\alpha}(\partial\Omega^i)$, and thus Lemma \ref{regularity} ensures that $\eta\in C^{0}(\partial\Omega^i)$. 

Similarly, to prove that $J_\lambda$ is invertible from $C^{0,\alpha}(\partial\Omega^i)$ to itself we observe that $J_\lambda$ is continuous from $C^{0,\alpha}(\partial\Omega^i)$ to itself and that $J_\lambda[\eta]\in C^{0,\alpha}(\partial\Omega^i)$ implies $\eta\in C^{0,\alpha}(\partial\Omega^i)$ for all $\eta\in L^2(\partial\Omega^i)$.\qed

\vspace{\baselineskip}

We now turn to consider problem \eqref{nlprob.small} for $\epsilon=0$. By the previous Lemma \ref{representation.small} and by the jump properties of the single and double layer potentials (cf.~equality \eqref{jump}) we deduce the following.

\begin{lem}\label{solution.small}
If $(\mu^o,\eta)\in C^{1,\alpha}(\partial\Omega^o)\times C^{0,\alpha}(\partial\Omega^i)$, then  the pair $(U^o[\mu^o,\eta,\eta], U^i[\mu^o,\eta,\eta])$ is a solution of \eqref{nlprob.small} with $\epsilon=0$ if and only if 
\begin{equation}\label{solution.small.eq1}
\left\{
\begin{array}{l}
(\frac{1}{2}I_{\Omega^o}+W_{\Omega^o})\mu^o=f^o-v^-_{\Omega^i}[\eta]_{|\partial\Omega^o}\,,\\
(\frac{1}{2}I_{\Omega^o}+\frac{\lambda-1}{\lambda+1}W^*_{\Omega^o})\eta+\frac{\lambda-1}{\lambda+1}\nu_{\Omega^i}\cdot\nabla w^+_{\Omega^o}[\mu^o]_{|\partial\Omega^i}\\
\qquad\qquad\qquad =\frac{\lambda}{\lambda+1}\mathcal{F}_G(\lambda^{-1}w^+_{\Omega^o}[\mu^o]_{|\partial\Omega^i}+\lambda^{-1} v^+_{\Omega^i}[\eta]_{|\partial\Omega^i})\,.
\end{array}
\right.
\end{equation}
\end{lem}

We show the existence of a solution of \eqref{solution.small.eq1} by an argument based on the invariance of the Leray-Schauder  topological degree (cf.~Theorem  \ref{topological}).

\begin{prop}\label{small} Assume that   $\mathcal{F}_G$ maps $C^{0,\alpha}(\partial\Omega^i)$ to itself and that  $|G(x,t)|\le C(1+|t|)^\delta$ for some  $C>0$, $\delta\in[0,1[$ and for all $(x,t)\in\partial\Omega^i\times\mathbb{R}$. Then there exists at least a solution $(\mu^o_0,\eta_0)\in C^{1,\alpha}(\partial\Omega^o)\times C^{0,\alpha}(\partial\Omega^i)$  of \eqref{solution.small.eq1}.
\end{prop}
\proof Since $\frac{1}{2}I_{\Omega^o}+W_{\Omega^o}$ is an invertible operator from $C^{1,\alpha}(\partial\Omega^o)$ to itself (cf.~Lemmas \ref{fredC} and \ref{ker}) it is enough to show that there exists a solution $\eta_0\in C^{0,\alpha}(\partial\Omega^i)$ of 
\begin{equation}\label{solution.small.eq2}
\begin{split}
&\eta=J_\lambda^{(-1)}\Bigg[\frac{\lambda}{\lambda+1}\mathcal{F}_G\Bigg(\frac{1}{\lambda}w^+_{\Omega^o}\Big[(\frac{1}{2}I_{\Omega^o}+W_{\Omega^o})^{(-1)}(f^o-v^-_{\Omega^i}[\eta]_{|\partial\Omega^o})\Big]_{|\partial\Omega^i}+\frac{1}{\lambda}v^+_{\Omega^i}[\eta]_{|\partial\Omega^i}\Bigg)\\
&\qquad\qquad\qquad\qquad-\frac{\lambda-1}{\lambda+1}\nu_{\Omega^i}\cdot\nabla w^+_{\Omega^o}\Big[(\frac{1}{2}I_{\Omega^o}+W_{\Omega^o})^{(-1)}f^o\Big]_{|\partial\Omega^i}
\Bigg]
\end{split}
\end{equation}  
(cf.~Lemma \ref{Jl}). We first show that the equation \eqref{solution.small.eq2} has a solution in $C^{0}(\partial\Omega^i)$. By the properties of integral operators with real analytic kernel and no singularities, by the invertibility of $\frac{1}{2}I_{\Omega^o}+W_{\Omega^o}$ in $C^0(\partial\Omega^o)$ (cf.~Lemmas \ref{fredC} and \ref{ker}), and by the mapping properties of the single layer potential  (cf., {\it e.g.},  Kress \cite[Thm.~2.22]{Kr89}, see also Miranda \cite[Chap.~II, \S14, III]{Mi70}) we verify that the map from $C^{0}(\partial\Omega^i)$ to itself which takes a function $\eta$ to 
\[
\frac{1}{\lambda}w^+_{\Omega^o}\left[(\frac{1}{2}I_{\Omega^o}+W_{\Omega^o})^{(-1)}(f^o-v^-_{\Omega^i}[\eta]_{|\partial\Omega^o})\right]_{|\partial\Omega^i}+\frac{1}{\lambda}v^+_{\Omega^i}[\eta]_{|\partial\Omega^i}
\]
is compact (cf.~Section \ref{sec.preliminaries}). In addition,   $\mathcal{F}_G$ is continuous from $C^{0}(\partial\Omega^i)$ to itself (because $G$ is continuous). It follows that the map which takes $\eta$ to 
\[
\mathcal{F}_G\left(\frac{1}{\lambda}w^+_{\Omega^o}\left[(\frac{1}{2}I_{\Omega^o}+W_{\Omega^o})^{(-1)}(f^o-v^-_{\Omega^i}[\eta]_{|\partial\Omega^o})\right]_{|\partial\Omega^i}+\frac{1}{\lambda}v^+_{\Omega^i}[\eta]_{|\partial\Omega^i}\right)
\]
is continuous from $C^{0}(\partial\Omega^i)$ to itself and maps bounded sets to sets with compact closure. Then, Lemma \ref{Jl} implies that the map from $C^{0}(\partial\Omega^i)$ to itself which takes $\eta$ to the right hand side of equation \eqref{solution.small.eq2} is continuous  and maps bounded sets to sets with compact closure. 
 Now let $t\in[0,1]$ and assume that 
\[
\begin{split}
&\eta=tJ_\lambda^{(-1)}\Bigg[\frac{\lambda}{\lambda+1}\mathcal{F}_G\Bigg(\frac{1}{\lambda}w^+_{\Omega^o}\Big[(\frac{1}{2}I_{\Omega^o}+W_{\Omega^o})^{(-1)}(f^o-v^-_{\Omega^i}[\eta]_{|\partial\Omega^o})\Big]_{|\partial\Omega^i}+\frac{1}{\lambda}v^+_{\Omega^i}[\eta]_{|\partial\Omega^i}\Bigg)\\
&\qquad\qquad\qquad\qquad-\frac{\lambda-1}{\lambda+1}\nu_{\Omega^i}\cdot\nabla w^+_{\Omega^o}\Big[(\frac{1}{2}I_{\Omega^o}+W_{\Omega^o})^{(-1)}f^o\Big]_{|\partial\Omega^i}
\Bigg]
\end{split}
\]
Then, by exploiting inequality $|G(x,t)|\le C(1+|t|)^\delta$ one verifies that
\begin{equation}\label{solution.small.eq3}
\|\eta\|_{C^{0}(\partial\Omega^i)} \le c_1\left(c_2+c_3\|\eta\|_{C^{0}(\partial\Omega^i)}\right)^\delta+c_4
\end{equation}
where $c_1$, \dots, $c_4$ are real positive numbers which depend on $C$, $t$, and $\lambda$, on the norm of the bounded operator $J^{(-1)}_\lambda$ from $C^{0}(\partial\Omega^i)$ to itself, on the norm of the bounded operator from $C^{0}(\partial\Omega^i)$ to itself which takes $\phi$ to 
\[
- w^+_{\Omega^o}\left[(\frac{1}{2}I_{\Omega^o}+W_{\Omega^o})^{(-1)}v^-_{\Omega^i}[\phi]_{|\partial\Omega^o}\right]_{|\partial\Omega^i}+v^+_{\Omega^i}[\phi]_{|\partial\Omega^i}\,,
\] and on the $C^0(\partial\Omega^i)$ norms of the functions $w^+_{\Omega^o}\left[(\frac{1}{2}I_{\Omega^o}+W_{\Omega^o})^{(-1)}f^o \right]_{|\partial\Omega^i}$ and 
$\nu_{\Omega^i}\cdot\nabla w^+_{\Omega^o}[(\frac{1}{2}I_{\Omega^o}+W_{\Omega^o})^{(-1)}f^o]_{|\partial\Omega^i}$.
Then inequality \eqref{solution.small.eq3} implies that 
\[
\|\eta\|_{C^{0}(\partial\Omega^i)} \le\max\left\{1,(c_1(c_2+c_3)^\delta+c_4)^{1/(1-\delta)}\right\}
\]
Thus Theorem \ref{topological} implies that there exists $\eta_0\in C^{0}(\partial\Omega^i)$ solution of \eqref{solution.small.eq2}. Then, by  classical results of potential theory (cf.~Miranda \cite[Chap.~II, \S14, III]{Mi70}), we have $v^+_{\Omega^i}[\eta_0]_{|\partial\Omega^i}\in C^{0,\alpha}(\partial\Omega^i)$, and, by the properties of integral operators with real analytic kernels and no singularities, by the assumption that $\mathcal{F}_G$ maps $C^{0,\alpha}(\partial\Omega^i)$ to itself, and by equation \eqref{solution.small.eq2} we deduce that $\eta_0\in C^{0,\alpha}(\partial\Omega^i)$. \qed

\vspace{\baselineskip}

We now pass to consider $\epsilon\neq 0$. We assume that 
\begin{equation}\label{frechet}
\begin{split}
&\text{the composition operators $\mathcal{F}_\Phi$ and $\mathcal{F}_G$  are continuously Fr\'echet differentiable }\\
&\text{from $C^{1,\alpha}(\partial\Omega^i)$ to itself and from $C^{0,\alpha}(\partial\Omega^i)$ to itself, respectively.}
\end{split}
\end{equation} 
We observe that condition \eqref{frechet} implies that the partial derivatives $\partial_t\Phi(x,t)$ and $\partial_t G(x,t)$ exist for all $(x,t)\in\partial\Omega^i\times\mathbb{R}$, that the composition operators $\mathcal{F}_{\partial_t\Phi}$ and $\mathcal{F}_{\partial_tG}$ map $C^{1,\alpha}(\partial\Omega^i)$ to itself and $C^{0,\alpha}(\partial\Omega^i)$ to itself, respectively,  and that
\[
\begin{split}
& d\mathcal{F}_\Phi(v_0).v= (\mathcal{F}_{\partial_t\Phi}\, v_0)\,v\qquad\forall v\in C^{1,\alpha}(\partial\Omega^i)\,,\\
& d\mathcal{F}_G(w_0).w=(\mathcal{F}_{\partial_tG}\, w_0)\,w\qquad\forall w\in C^{0,\alpha}(\partial\Omega^i)\,,
\end{split}
\] 
where $d\mathcal{F}_\Phi(v_0)$ denotes the differential of $\mathcal{F}_{\Phi}$ evaluated at a function $v_0\in C^{1,\alpha}(\partial\Omega^i)$ and $d\mathcal{F}_G(w_0)$ denotes  the differential of $\mathcal{F}_{G}$ evaluated at a function $w_0\in C^{0,\alpha}(\partial\Omega^i)$  (cf., {\it e.g.}, Lanza de Cristoforis \cite[Prop.~6.3]{La07}). 

Now we introduce the nonlinear operator $N\equiv(N^o,N^i_1,N^i_2)$ from $\mathbb{R}\times C^{1,\alpha}(\partial\Omega^o)\times C^{0,\alpha}(\partial\Omega^i)^2$ to $C^{1,\alpha}(\partial\Omega^o)\times C^{1,\alpha}(\partial\Omega^i)\times C^{0,\alpha}(\partial\Omega^i)$ which takes $(\epsilon,\mu^o,\eta^o,\eta^i)$ to 
\[
\begin{split}
&N^o[\epsilon,\mu^o,\eta^o,\eta^i]\equiv (\frac{1}{2}I_{\Omega^o}+W_{\Omega^o})\mu^o+v_{\Omega^i}^-[\eta^o]_{|\partial\Omega^o}-f^o\,,\\
&N^i_1[\epsilon,\mu^o,\eta^o,\eta^i]\equiv v_{\Omega^i}[\eta^o-\eta^i]_{|\partial\Omega^i}-\epsilon\mathcal{F}_\Phi(\lambda^{-1}w^+_{\Omega^o}[\mu^o]_{|\partial\Omega^i}+\lambda^{-1}v^+_{\Omega^i}[\eta^i]_{|\partial\Omega^i})\,,\\
&N^i_2[\epsilon,\mu^o,\eta^o,\eta^i]\equiv (\frac{1}{2}I_{\Omega^i}+ W^*_{\Omega^i})\eta^o-\lambda^{-1}(-\frac{1}{2}I_{\Omega^i}+ W^*_{\Omega^i})\eta^i\\
&\qquad\qquad\qquad\quad +(\lambda-1)\lambda^{-1} \nu_{\Omega^i}\cdot\nabla w^+_{\Omega^o}[\mu^o]_{|\partial\Omega^i}-\mathcal{F}_G(\lambda^{-1}w^+_{\Omega^o}[\mu^o]_{|\partial\Omega^i}+\lambda^{-1}v^+_{\Omega^i}[\eta^i]_{|\partial\Omega^i})\,.
\end{split}
\]
Then, by the mapping and jump properties of single and double layer potentials (cf.~Section \ref{sec.preliminaries}), one  verifies the validity of the following Lemmas \ref{Nanalytic} and \ref{Nsolution}.

\begin{lem}\label{Nanalytic}
 If $\Phi$ and $G$ satisfy condition  \eqref{frechet}, then $N$ is continuously Fr\'echet differentiable map  from $\mathbb{R}\times C^{1,\alpha}(\partial\Omega^o)\times  C^{0,\alpha}(\partial\Omega^i)^2$ to $C^{1,\alpha}(\partial\Omega^i)\times C^{0,\alpha}(\partial\Omega^i)$.
\end{lem}

\begin{lem}\label{Nsolution}
Let $(\epsilon,\mu^o,\eta^o,\eta^i)\in \mathbb{R}\times C^{1,\alpha}(\partial\Omega^o)\times C^{0,\alpha}(\partial\Omega^i)^2$. Then $N[\epsilon,\mu^o,\eta^o,\eta^i]=0$ if and only if $(U^o[\mu^o,\eta^o,\eta^i], U^i[\mu^o,\eta^o,\eta^i])$ is a solution of \eqref{nlprob.small}.
\end{lem}

Moreover, one can prove the following.

\begin{lem}\label{frechet.lem}
Let  $\Phi$ and $G$ satisfy condition  \eqref{frechet}.  Let $(\tilde\mu^o,\tilde\eta^o,\tilde\eta^i)\in C^{1,\alpha}(\partial\Omega^i)\times C^{0,\alpha}(\partial\Omega^i)^2$. If 
\begin{equation}\label{ass.G}
{\partial_t G}(x,\lambda^{-1}w^+_{\Omega^o}[\tilde\mu^o](x)+\lambda^{-1}v^+_{\Omega^i}[\tilde\eta^i](x))\ge 0\qquad\forall x\in\partial\Omega^i,
\end{equation} 
then $\partial_{(\mu^o,\eta^o,\eta^i)}N[0,\tilde\mu^o,\tilde\eta^o,\tilde\eta^i]$ (the partial differential  of $N$  with respect to $(\mu^o,\eta^o,\eta^i)$ evaluated at $(0,\tilde\mu^o,\tilde\eta^o,\tilde\eta^i)$) is an isomorphism from $C^{1,\alpha}(\partial\Omega^o)\times C^{0,\alpha}(\partial\Omega^i)^2$ to $C^{1,\alpha}(\partial\Omega^o)\times C^{1,\alpha}(\partial\Omega^i) \times C^{0,\alpha}(\partial\Omega^i)$.
\end{lem}
\proof We have 
\[
\begin{split}
&\partial_{(\mu^o,\eta^o,\eta^i)}N^o[0,\tilde\mu^o,\tilde\eta^o,\tilde\eta^i](\bar\mu^o,\bar\eta^o,\bar\eta^i)=(\frac{1}{2}I_{\Omega^o}+W_{\Omega^o})\bar\mu^o+v_{\Omega^i}^-[\bar\eta^o]_{|\partial\Omega^o}\,,\\
&\partial_{(\mu^o,\eta^o,\eta^i)}N^i_1[0,\tilde\mu^o,\tilde\eta^o,\tilde\eta^i](\bar\mu^o,\bar\eta^o,\bar\eta^i)= v_{\Omega^i}[\bar\eta^o-\bar\eta^i]_{|\partial\Omega^i}\,,\\
&\partial_{(\mu^o,\eta^o,\eta^i)}N^i_2[0,\tilde\mu^o,\tilde\eta^o,\tilde\eta^i](\bar\mu^o,\bar\eta^o,\bar\eta^i)= (\frac{1}{2}I_{\Omega^i}+ W^*_{\Omega^i})\bar\eta^o-\lambda^{-1}(-\frac{1}{2}I_{\Omega^i}+ W^*_{\Omega^i})\bar\eta^i\\
&\qquad\qquad\qquad\qquad +(\lambda-1)\lambda^{-1} \nu_{\Omega^i}\cdot\nabla w^+_{\Omega^o}[\bar\mu^o]_{|\partial\Omega^i} -\tilde\gamma\,(\lambda^{-1}w^+_{\Omega^o}[\bar\mu^o]_{|\partial\Omega^i}+\lambda^{-1}v^+_{\Omega^i}[\bar\eta^i]_{|\partial\Omega^i})\,,
\end{split}
\] 
for all $(\bar\mu^o,\bar\eta^o,\bar\eta^i)\in C^{1,\alpha}(\partial\Omega^i)\times C^{0,\alpha}(\partial\Omega^i)^2$ where  $\tilde\gamma$ is the function of $C^{0,\alpha}(\partial\Omega^i)$  defined by
\[
\tilde\gamma(x)\equiv{\partial_t G}(x,w^+_{\Omega^o}[\tilde\mu^o](x)+v^-_{\Omega^i}[\tilde\eta^i](x))\qquad\forall x\in\partial\Omega^i.
\]
Then we observe that the operator which takes $(\bar\mu^o,\bar\eta^o,\bar\eta^i)$ to $\partial_{(\mu^o,\eta^o,\eta^i)}N[0,\tilde\mu^o,\tilde\eta^o,\tilde\eta^i](\bar\mu^o,\bar\eta^o,\bar\eta^i)$ is Fredholm of index $0$ from $C^{1,\alpha}(\partial\Omega^i)\times C^{0,\alpha}(\partial\Omega^i)^2$ to itself. Indeed the operator which takes $(\bar\mu^o,\bar\eta^o,\bar\eta^i)$ to 
\[
\Big((\frac{1}{2}I_{\Omega^o}+W_{\Omega^o})\bar\mu^o\,,\; v_{\Omega^i}[\bar\eta^o-\bar\eta^i]_{|\partial\Omega^i}\,,\; \frac{1}{2}(\bar\eta^o+\lambda^{-1}\bar\eta^i)\Big)
\]
is an isomorphism (cf.~Lemmas \ref{fredC} and \ref{ker} and condition \eqref{ass.omega}) and the operator which takes $(\bar\mu^o,\bar\eta^o,\bar\eta^i)$ to 
\[
\begin{split}
&\Big(v_{\Omega^i}^-[\bar\eta^o]_{|\partial\Omega^o}\,,\; 0\,,\; W^*_{\Omega^i}[\bar\eta^o]-\lambda^{-1}W^*_{\Omega^i}[\bar\eta^i]\\
&\qquad\qquad\qquad\quad +(\lambda-1)\lambda^{-1} \nu_{\Omega^i}\cdot\nabla w^+_{\Omega^o}[\bar\mu^o]_{|\partial\Omega^i} - \tilde\gamma(\lambda^{-1}w^+_{\Omega^o}[\bar\mu^o]_{|\partial\Omega^i}+\lambda^{-1}v^+_{\Omega^i}[\bar\eta^i]_{|\partial\Omega^i})\Big)
\end{split}
\] 
is compact (by the properties of integral operators with real analytic kernels and no singularity and by Lemma \ref{schauder}). Hence, to prove the statement of the lemma, it suffices to show that 
\[
\partial_{(\mu^o,\eta^o,\eta^i)}N[0,\tilde\mu^o,\tilde\eta^o,\tilde\eta^i](\bar\mu^o,\bar\eta^o,\bar\eta^i)=0
\]
implies $(\bar\mu^o,\bar\eta^o,\bar\eta^i)=(0,0,0)$. If $\partial_{(\mu^o,\eta^o,\eta^i)}N[0,\tilde\mu^o,\tilde\eta^o,\tilde\eta^i](\bar\mu^o,\bar\eta^o,\bar\eta^i)=0$, then by the jump properties of the single and double layer potentials the pair  $(U^o[\bar\mu^o,\bar\eta^o,\bar\eta^i], U^i[\bar\mu^o,\bar\eta^o,\bar\eta^i])$ is a solution of the problem
\[
\left\{
\begin{array}{ll}
\Delta u^o=0&\text{in }\Omega^o\setminus\mathrm{cl}\Omega^i\,,\\
\Delta u^i=0&\text{in }\Omega^i\,,\\
u^o(x)=0&\text{for all }x\in\partial\Omega^o\,,\\
u^o(x)=\lambda u^i(x)&\text{for all }x\in\partial\Omega^i\,,\\
\nu_{\Omega^i}\cdot\nabla u^o(x)-\nu_{\Omega^i}\cdot\nabla u^i(x)=\tilde\gamma(x)
\, u^i(x)&\text{for all }x\in\partial\Omega^i\,.\\
\end{array}
\right.
\] 
Then, by inequalities $\lambda> 0$ and $\tilde\gamma\ge 0$ (cf.~condition \eqref{ass.G}) and by a standard energy argument one verifies that $(U^o[\bar\mu^o,\bar\eta^o,\bar\eta^i], U^i[\bar\mu^o,\bar\eta^o,\bar\eta^i])=(0,0)$. Thus $(\bar\mu^o,\bar\eta^o,\bar\eta^i)=(0,0,0)$ by Lemma \ref{representation.small} and the proof is completed.
\qed

\vspace{\baselineskip}

Then, by Lemma \ref{frechet.lem} and by the implicit function theorem (see, {\it e.g.}, Deimling \cite[\S15]{De85}) one verifies the validity of the following proposition.

\begin{prop}\label{implicit.prop}
Let $\Phi$ and $G$ satisfy \eqref{frechet}. Let $(\tilde\mu^o,\tilde\eta^o,\tilde\eta^i)\in C^{1,\alpha}(\partial\Omega^i)\times C^{0,\alpha}(\partial\Omega^i)^2$ and $N[0, \tilde\mu^o,\tilde\eta^o,\tilde\eta^i]=0$. Assume that condition \eqref{ass.G} holds true. Then there exist $\epsilon^*>0$, a neighbourhood $\mathcal{U}$ of $(\tilde\mu^o,\tilde\eta^o,\tilde\eta^i)$ in  $C^{1,\alpha}(\partial\Omega^i)\times C^{0,\alpha}(\partial\Omega^i)^2$, and  a continuously Fr\'echet differentiable map $(\mu^o[\cdot],\eta^o[\cdot],\eta^i[\cdot])$ from $]-\epsilon^*,\epsilon^*[$ to $\mathcal{U}$ such that the set of zeros of $N$ in $]-\epsilon^*,\epsilon^*[\times\mathcal{U}$ coincides with the graph of $(\mu^o[\cdot],\eta^o[\cdot],\eta^i[\cdot])$. In particular, $N[\epsilon, \mu^o[\epsilon],\eta^o[\epsilon],\eta^i[\epsilon]]=0$ for all $\epsilon\in]-\epsilon^*,\epsilon^*[$ and $(\mu^o[0],\eta^o[0],\eta^i[0])=(\tilde\mu^o,\tilde\eta^o,\tilde\eta^i)$.
\end{prop}

We are now ready to prove the main Theorem \ref{thm.small} of this section.

\begin{thm}\label{thm.small} Let $\Phi$ and $G$ satisfy condition \eqref{frechet}. Assume that $|G(x,t)|\le C(1+|t|)^\delta$ for some  $C>0$, $\delta\in[0,1[$ and for all $(x,t)\in\partial\Omega^i\times\mathbb{R}$. Then the following statement hold:
\begin{enumerate}
\item[(i)] there exists at least a solution $(u^o_0,u^i_0)\in C^{1,\alpha}(\mathrm{cl}\Omega^o\setminus\Omega^i)\times C^{1,\alpha}(\mathrm{cl}\Omega^i)$ of the boundary value problem in \eqref{nlprob.small} with $\epsilon=0$. 
\end{enumerate}
If in addition we have 
\begin{equation}\label{ass.dGu}
(\partial_ tG)(x,u^i_0(x))\ge 0\qquad\forall x\in\partial\Omega^i\,
\end{equation}
then there exist $\epsilon_*>0$ and a family of functions $\{(u^o_\epsilon,u^i_\epsilon)\}_{\epsilon\in]-\epsilon_*,\epsilon_*[\setminus\{0\}}$ such that following statements hold:
\begin{enumerate}
\item[(ii)]  for all $\epsilon\in]-\epsilon_*,\epsilon_*[$ the pair $(u^o_\epsilon,u^i_\epsilon)$ belongs to $C^{1,\alpha}(\mathrm{cl}\Omega^o\setminus\Omega^i)\times C^{1,\alpha}(\mathrm{cl}\Omega^i)$ and is a solution of \eqref{nlprob.small};
\item[(iii)] the map from $]-\epsilon_*,\epsilon_*[$ to $C^{1,\alpha} (\mathrm{cl}\Omega^o\setminus\Omega^i)\times C^{1,\alpha}(\mathrm{cl}\Omega^i)$ which takes $\epsilon$ to $(u^o_\epsilon,u^i_\epsilon)$ is continuously Fr\'echet differentiable;
\item[(iv)] there exists an open subset $\mathcal{V}$ of $C^{1,\alpha}_\mathrm{harm}(\mathrm{cl}\Omega^o\setminus\Omega^i)\times C^{1,\alpha}_\mathrm{harm}(\mathrm{cl}\Omega^i)$ such that, for all  fixed $\epsilon\in]-\epsilon_*,\epsilon_*[$ the pair $(u^o_\epsilon,u^i_\epsilon)$ is the unique solution  of \eqref{nlprob.small} belonging to $\mathcal{V}$.
\end{enumerate}
\end{thm}
\proof  (i) By Lemma \ref{small} there exists at least a solution $(\mu^o_0,\eta_0)\in C^{1,\alpha}(\partial\Omega^o)\times C^{0,\alpha}(\partial\Omega^i)$ of \eqref{solution.small.eq1}. Then we define $(u^o_0,u^i_0)\equiv(U^o[\mu^o_0,\eta_0,\eta_0],U^i[\mu^o_0,\eta_0,\eta_0])$  and the validity of statement (i) follows by Lemma \ref{solution.small} (see also Proposition \ref{implicit.prop}). 

(ii) Since $(\mu^o_0,\eta_0)$ is a solution  of \eqref{solution.small.eq1}, we have $N[0,\mu^o_0,\eta_0,\eta_0]=0$. Then let $(\tilde\mu^o,\tilde\eta^o,\tilde\eta^i)\equiv(\mu^o_0,\eta_0,\eta_0)$.  By condition \eqref{ass.dGu} and by the jump properties of single and double layer potential, one verifies that condition \eqref{ass.G} is satisfied. Accordingly, the assumption of Proposition \ref{implicit.prop} are fulfilled and we can take $\epsilon_*\equiv \epsilon^*$ and define $(u^o_\epsilon,u^i_\epsilon)\equiv(U^o[
\mu^o[\epsilon],\eta^o[\epsilon],\eta^i[\epsilon]],U^i[\mu^o[\epsilon],\eta^o[\epsilon],\eta^i[\epsilon]])$ for all $\epsilon\in]-\epsilon_*,\epsilon_*[$. The validity of (ii) follows by Lemma \ref{Nsolution}.

(iii)  It is a consequence of the continuous Fr\'echet differentiability of  $(\mu^o[\cdot],\eta^o[\cdot],\eta^i[\cdot])$, of the definition of $(U^o,U^i)$ in Lemma \ref{representation.small}, and of the mapping properties of the single and double layer potentials (cf.~Miranda \cite{Mi65}).

(iv) Let $\mathcal{U}$ be the open neighbourhood of $(\tilde\mu^o,\tilde\eta^o,\tilde\eta^i)$ introduced in Proposition \ref{implicit.prop}. Let $\mathcal{V}\equiv \{(U^o[\mu^o,\eta^o,\eta^i]\,,\;U^i[\mu^o,\eta^o,\eta^i])\,:\;(\mu^o,\eta^o,\eta^i)\in\mathcal{U}\}$. Since $U=(U^o,U^i)$ is an open operator the set $\mathcal{V}$ is open in $C^{1,\alpha}_\mathrm{harm}(\mathrm{cl}\Omega^o\setminus\Omega^i)\times C^{1,\alpha}_\mathrm{harm}(\mathrm{cl}\Omega^i)$ (cf.~Lemma \ref{representation.small}). Moreover, the pair $(u^o_\epsilon,u^i_\epsilon)=(U^o[
\mu^o[\epsilon],\eta^o[\epsilon],\eta^i[\epsilon]],U^i[\mu^o[\epsilon],\eta^o[\epsilon],\eta^i[\epsilon]])$ belongs to $\mathcal{V}$ for all $\epsilon\in]-\epsilon_*,\epsilon_*[$ (cf.~Proposition \ref{implicit.prop}). Now fix $\epsilon_\sharp\in]-\epsilon_*,\epsilon_*[$  and assume that $(u^o_\sharp,u^i_\sharp)\in\mathcal{V}$ is a solution of \eqref{nlprob.small} for $\epsilon=\epsilon_\sharp$. Then there exists $(\mu^o_\sharp,\eta^o_\sharp,\eta^i_\sharp)\in \mathcal{U}$ such that  $(u^o_\sharp,u^i_\sharp)=(U^o[\mu^o_\sharp,\eta^o_\sharp,\eta^i_\sharp],U^i[\mu^o_\sharp,\eta^o_\sharp,\eta^i_\sharp])$. Moreover, $N[\epsilon_\sharp,\mu^o_\sharp,\eta^o_\sharp,\eta^i_\sharp]=0$ by Lemma \ref{Nsolution} and thus $(\mu^o_\sharp,\eta^o_\sharp,\eta^i_\sharp)=(\mu^o[\epsilon_\sharp],\eta^o[\epsilon_\sharp],\eta^i[\epsilon_\sharp])$ by Proposition \ref{implicit.prop}. Accordingly $(u^o_\sharp,u^i_\sharp)=(u^o_{\epsilon_\sharp},u^i_{\epsilon_\sharp})$ and the proof is complete.\qed

\vspace{\baselineskip}

\section*{Acknowledgment}
The research  of M.~Dalla Riva was supported by the Portuguese funds through the CIDMA - Center for Research and Development in Mathematics and Applications, and the Portuguese Foundation for Science and Technology (``FCT--Funda{\c c}{\~a}o para a Ci\^encia e a Tecnologia''), within project PEst-OE/MAT/UI4106/2014. The research of M.~Dalla Riva was also supported by the Portuguese Foundation for Science and Technology (``FCT--Funda{\c c}{\~a}o para a Ci\^encia e a Tecnologia'')    with the research grant SFRH/BPD/ 64437/2009. In addition, the work of M.~Dalla Riva was supported by ``Progetto di Ateneo: Singular perturbation problems for differential operators -- CPDA120171/12'' of the University of Padova. 

G.~Mishuris acknowledges the support of the European Community's Seven Framework Programme under contract number PIAPP-GA-284544-PARM-2.

Finally, M. Dalla Riva wishes to thank the Department of Mathematics of the University of Aberystwyth for the hospitality received during the development of a part of the work.

\end{document}